\documentclass[francais]{smfart}
\usepackage{smfthm}
\theoremstyle{plain}

\def\k#1{\kern#1em}
\def\Ib#1{{I\kern-.25em#1}}
\def\Ibb#1{{I\kern-.23em#1}}
\def\vh#1{\vrule  width.02em height#1ex depth0ex}
\def\vb#1{\vrule  width.02em height1.47ex depth#1ex}
\def\vcg{\vrule  width.02em height1.4ex depth-.05ex}
\def\AA{{\rm{\k{.22}\vh{.5}\k{-.24}A}}}

\def\CC{{\rm{\k{.24}\vcg\k{-.26}C}}}

\def\NN{{\rm\Ibb N}}

\def\RR{{\rm\Ib R}}

\def\ZZ{{\rm{\k{.26}\vh{0.5}\k{.04}\vb{-1}\k{-.34}Z}}}
\def\ind{{\rm{ 1 \hspace {-.1cm} I}}}

\author{Fr\'ed\'eric Sarkis}
\address{Institut de math\'ematiques \\
Universit\'e Pierre et Marie Curie, Tour 46-56 B. 507 \\
4 Place Jussieu, 75005 Paris}
\email{sarkis@math.jussieu.fr}
\title{Probl\`eme Plateau complexe dans les vari\'et\'es K\"ahl\'eriennes}
\alttitle{Complex Plateau problem in K\"ahler manifolds}

\begin{document}

\begin{abstract}
L'\'etude du ``probl\`eme Plateau complexe'' (ou ``probl\`eme du bord'') dans une vari\'et\'e complexe $X$
consiste \`a caract\'eriser les sous-vari\'et\'es r\'eelles $\Gamma$ de
$X$ qui sont le bord de sous-ensembles analytiques de 
$\omega\backslash \Gamma$. Notre principal r\'esultat traite le cas 
$X=U \times \omega$ o\`u $U$ est une vari\'et\'e complexe connexe et
$\omega$ est une vari\'et\'e K\"ahl\'erienne disque convexe.
Comme cons\'equence, nous obtenons des 
r\'esultats de Harvey-Lawson \cite{harvey2}, Dolbeault-Henkin
\cite{henkin} et Dinh \cite{dinh2}. Nous obtenons aussi
une  g\'en\'eralisation des th\'eor\`emes de Hartogs-Levi
et Hartogs-Bochner.
Finalement, nous montrons qu'une structure CR strictement
pseudo-convexe plongeable dans une vari\'et\'e K\"ahl\'erienne
disque-convexe est plongeable dans $\CC^n$ si
et seulement si elle admet une fonction CR non constante.
\end{abstract}

\begin{altabstract}
The ``complex Plateau problem'' (or ``boundary problem'') in a complex
manifold $X$ is the problem of characterizing the real submanifolds
$\Gamma$ of $X$ which are boundaries of analytic sub-varieties of $X
\backslash \Gamma$. Our principal result treat the case $X=U\times
\omega$ where $U$ is a connected complex manifold and $\omega$ is a
disk-convex K\"ahler manifold. As a consequence, we obtain
 results of Harvey-Lawson \cite{harvey2}, Dolbeault-Henkin
\cite{henkin} and Dinh \cite{dinh2}.
We also give a generalization of
Hartogs-Levi and Hartogs-Bochner theorems.
Finally, we prove that a
strictly pseudo-convex CR structure embeddable in a disk-convex
K\"ahler manifold is embeddable in $\CC^n$ 
if and only if it has a non constant CR function. 
\end{altabstract}

\maketitle

\maketitle
\section { Introduction.}
L'\'etude du ``Probl\`eme Plateau complexe'' (ou ``probl\`eme du bord'') dans une vari\'et\'e complexe $X$
consiste \`a caract\'eriser les sous-vari\'et\'es r\'eelles $X$ de
$X$ qui sont le bord (au sens des courants) 
de sous-ensembles analytiques de $\omega
\backslash \Gamma$.\par

Dans l'espace affine, les courbes bords de surfaces de Riemann sont
caract\'eris\'ees par une condition int\'egrale nomm\'ee ``condition des
moments'' (voir \cite{wermer,stolzenberg,alexander,lawrence,dinh}).
Les sous-vari\'et\'es compactes, ferm\'ees et 
de dimension sup\'erieure ou \'egale
\`a trois qui sont bords d'ensembles analytiques de $\CC^n$ sont
celles dont la dimension de l'espace complexe tangent est maximal
en chaque point (voir \cite{harvey,chirka,dinh}).\par

Dans l'espace projectif, sous la condition pr\'ec\'edente, 
le probl\`eme du bord 
n'admet pas toujours de solution.
 Harvey et Lawson \cite{harvey2} 
ont cependant donn\'e une ca\-rac\-t\'erisation en terme de condition
des moments pour le probl\`eme du bord dans $P_n(\CC)\backslash
P_{n-p}(\CC)$ 
(o\`u $2p-1$ ($p \geq 2$) est la dimension de la
vari\'et\'e consid\'er\'ee). \par\noindent
R\'ecemment, Dolbeault et Henkin \cite{henkin} (puis Dinh \cite{dinh,dinh2})
ont donn\'e une condition n\'ecessaire et su\-ffi\-sante:
le probl\`eme du bord  pour une vari\'et\'e
maximalement complexe $M \subset P_n(\CC)$ admet une solution
s'il en admet une pour un nombre ``assez grand'' de tranches de $M$ 
par des sous-espaces lin\'eaires.\par
Le but de cet article est de g\'en\'eraliser
 ce dernier r\'esultat \`a une vari\'et\'e produit $U \times \omega$ 
o\`u $U$ est une vari\'et\'e complexe connexe et $\omega$ est une
vari\'et\'e  K\"ahl\'erienne compacte ou plus g\'en\'eralement {\it disque convexe} 
(i.e. pour tout compact $K \subset X$, 
il existe un compact $\widehat{K} \subset X$ 
tel que pour toute surface de Riemann $S$ et pour toute
application m\'eromorphe $f: S \rightarrow X$ telle que
$f(\partial S) \subset K$ on ait $f(S) \subset \widehat{K}$).
L'\'etude du probl\`eme du bord dans les espaces produits n'est pas
restrictive. En effet, le probl\`eme du bord dans l'espace projectif
peut toujours \^etre r\'eduit \`a l'\'etude du probl\`eme du bord dans
un espace produit. Comme nous le verrons, la r\'eciproque n'est
en g\'en\'eral pas vraie. 
De plus, cette r\'esolution parait la plus
adapt\'ee pour l'\'etude de l'extension des applications CR car 
si $f$ est une application CR, le graphe de $f$ est
naturellement dans un espace produit. 
Soit $X$ une vari\'et\'e k\"ahl\'erienne disque convexe,
nous obtenons alors les
corollaires suivants:
\begin{enumerate}
\item Une nouvelle d\'emonstration de la caract\'erisation g\'eom\'etrique 
du probl\`eme du bord dans $P_n(\CC)$ 
donn\'ee dans \cite{henkin, dinh2}.
\item Une g\'en\'eralisation du th\'eor\`eme de Hartogs-Levi 
pour les applications m\'ero\-mor\-phes \`a valeurs dans $X$.
En particulier, nous retrouvons les g\'en\'eralisations donn\'ees 
dans \cite{dinh2} et \cite{Ivashkovich}.
\item Une g\'en\'eralisation du th\'eor\`eme de Hartogs-Bochner
pour les applications CR \`a valeurs dans les vari\'et\'es
K\"ahl\'eriennes disque convexes  (pour $X=P_n(\CC)$, nous retrouvons des r\'esultats de \cite{henkin,
porten, sarkis}).
\item La g\'en\'eralisation suivante du th\'eor\`eme de Hartogs-Bochner:
Supposons que $X$ est une vari\'et\'e K\"ahl\'erienne 
de dimension $2$ ne contenant aucune surface de Riemann compacte.
Soit $M$  une hypersurface
r\'eelle de $X$ la s\'eparant  en deux composantes connexes
$\Omega_1$ et $\Omega_2$. Alors  toute fonction holomorphe au
voisinage de $M$ admet une extension holomorphe sur $\Omega_1$ 
ou sur $\Omega_2$.
\end{enumerate}
Ainsi que le corollaire principal suivant nouveau m\^eme
dans le cas o\`u $X=P_n(\CC)$:

\begin{enumerate}
\item[4.] Soit $M$ une structure CR abstraite, strictement
pseudoconvexe  et plongeable dans $X$. 
La vari\'et\'e $M$ est plongeable dans $\CC^n$ (ou de mani\`ere
\'equivalente, admet une solution au probl\`eme du bord dans $X$) 
si et seulement si elle admet une fonction CR non constante.
\end{enumerate}
Ce dernier r\'esultat donne une nouvelle r\'eduction
du  probl\`eme suivant
(voir \cite{henkin}):
{\it Soit $M$ une vari\'et\'e strictement pseudoconvexe de $P_n(\CC)$.
La vari\'et\'e $M$ est-elle toujours le bord d'une $p$-cha\^\i ne
holomorphe (ou de mani\`ere \'equivalente $M$ est elle plongeable dans
l'espace affine) ? }  
En particulier, cela prouve de mani\`ere imm\'ediate que les exemples de 
structures CR d'Andreotti-Rossi \cite{rossi1,rossi2,andreotti,
  grauert,falbel} et de Barrett \cite{barrett} (cas o\`u l'espace des
fonctions CR est de dimension 1) ne sont plongeables dans aucune vari\'et\'e k\"ahl\'erienne disque
convexe $X$. En effet, ces vari\'et\'es ne sont pas plongeables dans l'espace
affine mais, par construction, elles admettent des fonctions CR non
constantes, on en d\'eduit alors directement 
 qu'elles ne sont pas plongeables dans $X$
(pour $X=P_n(\CC)$ et $M$ la structure CR d'Andreotti-Rossi, ceci est
d\'emontr\'e de mani\`ere diff\'erente dans \cite{sarkis}).



\section{Notations et r\'esultats pr\'eliminaires.}
\subsection{ Probl\`eme du bord non born\'e dans $\Delta_n \times
  \overline \Delta_m$.}

Notons $\Delta_q$ le polydisque unit\'e de $\CC^q$.
Nous nous int\'eressons ici au probl\`eme du bord pour les vari\'et\'es
maximalement complexe de dimension $2n+1$ de $\Delta_n \times \CC^m$ 
dont la projection sur $\CC^m$ est born\'ee.
En g\'en\'eral il n'existe pas de solution au probl\`eme du bord pour
de telles vari\'et\'es.
En effet, il suffit de prendre un courbe ferm\'ee 
$\gamma \subset \Delta_m$ qui n'est pas le bord d'une surface de
Riemann et de poser $\Gamma=\gamma \times \CC^n$, pour obtenir un contre-exemple.
Il est tout de m\^eme possible de donner une condition n\'ecessaire et
suffisante. En fait, cette condition est d\'ej\`a 
implicitement utilis\'ee dans \cite{dinh} et 
nous nous bornons ici \`a reformuler  
les propositions de \cite{dinh} dans ce cadre.

On dit qu'un compact $A\subset \RR^n$ {\it est g\'eom\'etriquement
  $m$-rectifiable} ou encore {\it est de classe $A_m$} si
$A$ est $({\mathcal H}_m,m)$-rectifiable et si le c\^one tangentiel de $A$ en
${\mathcal H}_m$-presque tout point est un espace vectoriel r\'eel de dimension
$m$ (o\`u  ${\mathcal H}_m$ est la mesure de Hausdorff de dimension $m$).\par

Soit $X$ un espace complexe, on appelle {\it $p$-cha\^{\i}ne
  holomorphe de $X$ } toute somme localement finie $[T]=\sum n_j [V_j]$
\`a coefficients $n_j$ dans $\ZZ$ de sous-ensembles analytiques $V_j$ 
de dimension $p$ de $X$.
 On appelle {\it volume de $[T]$}, l'expression
$$ \mbox{Vol }[T]=\sum |n_j|\mbox{Vol }V_j$$
o\`u $\mbox{Vol }V_j$ est le volume $2p$-dimensionnel de l'ensemble
analytique $V_j$; $\mbox{Vol }[T]$ est aussi la masse du courant
$[T]$. On notera $T$ le support de la $p$-cha\^{\i}ne holomorphe $[T]$
(i.e. $T=\cup_{\{j ; n_j\neq0\}}V_j$). Dans la suite, si $[T]$ est un
courant on notera aussi $Vol [T]$ la masse de ce courant.

Soient $X,Y$ deux vari\'et\'es r\'eelles lisses, $Y$ de dimension $p \leq m$ et
$f:X \rightarrow Y$ une application $C^\infty$ et $\Gamma$ un courant plat
de dimension $m$ (en particulier les courants rectifiables sont
plats).  Alors pour ${\mathcal H}_p$-presque tout $y \in Y$, 
la tranche $[\Gamma,f,y]$ est un
courant plat de dimension $m-p$ de support inclus dans $\Gamma \cap
f^{-1}(y)$ v\'erifiant:
$$\int_Y\Phi(y)([\Gamma,f,y],\Psi)d{\mathcal H}_p(y)=([\Gamma, f^*(\Phi \wedge
\Omega)]\wedge \Psi)$$ o\`u $\Psi$ est une $(m-p)$-forme \`a support compact,
$\Phi$ une fonction \`a support compact et $\Omega$ est la forme volume
de $Y$ (voir \cite{federer}).\par
Soit $[\Gamma]$ un
courant rectifiable dont le support $\Gamma$ est de classe $A_{2n+1}$ de
$\CC^{n+m}$. Alors,  pour presque tous les $m$-plan $\CC^n_\nu \subset
\CC^{n+m}$ et pour presque 
tous les points $z\in \CC^n_\nu$, 
la tranche $[\Gamma,\pi_\nu,z]$ (o\`u $pi_\nu$ est
la projection orthogonale sur $\CC^n_\nu$)
est un courant rectifiable dont le support
est de classe $A_{1}$ (voir \cite{dinh} Lemme 1.4). \par

Notons $\pi:\Delta_n \times \CC^m \rightarrow \Delta_n$ la projection
$(z_1,...,z_n)\times (w_1,...,w_m) \mapsto (z_1,...,z_n)$.
Soit $[\Gamma]$ un courant rectifiable, ferm\'e et maximalement
complexe de $\Delta_n \times \CC^m$ 
dont le support $\Gamma$ est de classe $A_{2n+1}$ et admet une
projection born\'ee sur $\CC^m$.  
Pour presque tout $z \in \Delta_n$, la tranche $[\Gamma, \pi, z]$ 
est un courant rectifiable ferm\'e de dimension 1. 

Soient $\chi_{\epsilon_j}$ des fonctions
de classe $C^\infty$, d\'efinies sur $\CC$, $\chi_{\epsilon_j}(x)=0$
  pour $|x| <\k{-.3}< \epsilon_j$ et $\chi_{\epsilon_j}=\frac{1}{2\pi i}$ pour
  $|x| > \epsilon$.

\begin{prop}\cite{dinh}
Pour tout $z \in \Delta_n$, pour toute $(1,0)$-forme $\phi$, holomorphe
 au voisinage de $\Gamma$, la fonction 
$${\mathcal{M}}(z,\phi)=\left([\Gamma],\phi\wedge\bigwedge_{j=1}^{n}
d\chi_{\epsilon_j}(\zeta_j-z_j)\wedge\frac {d\zeta_j}{ \zeta_j-z_j}\right)$$
est ind\'ependante des $\chi_{\epsilon_j}$ 
et \'egale \`a $([\Gamma, \pi, z],\phi)$ quand cette derni\`ere est
bien d\'efinie.
En particulier la fonction ${\mathcal{M}}(z,\phi)$ est holomorphe dans $\Delta_n$.
\end{prop}

Cette proposition  montre que la condition des moments (cas o\`u $\phi$
est une $(1,0)$-forme holomorphe) sur les tranches
d'une vari\'et\'e maximalement complexe varie
de mani\`ere holomorphe.
Le principe du prolongement
analytique nous dit alors que si la condition des moments est v\'erifi\'ee
pour un sous-ensemble assez grand $K \subset \Delta_n$, 
elle est v\'erifi\'ee sur tout $\Delta_n$.
Il est alors possible de recoller ces tranches
pour trouver une solution au  probl\`eme du bord 
pour $[\Gamma]$.

\begin{prop}
Soit $K$ un sous-ensemble de $\Delta_n$, non inclus dans une r\'eunion
localement finie d'ensembles analytiques de dimension $(n-1)$
immerg\'es dans $\Delta_n$.
Supposons que pour tout $\nu \in K$,
$[\Gamma, \pi, \nu]$ est bien d\'efinie et est un $1$-courant
rectifiable 
ferm\'e dont le support est de classe $A_1$
et que pour toute $2$-forme lisse $\Theta$ sur $\CC^m$ on aie
$[\Gamma]\wedge \Theta=0$. 
 Alors, les propositions
suivantes sont \'equivalentes:
\begin{enumerate}
\item  ${\mathcal M}(\nu,\phi)$ est nulle pour tout $\nu \in K$   
et toute $(1,0)$-forme holomorphe $\phi$ (i.e. pour tout $\nu \in K$,
$[\gamma_\nu]$ v\'erifie la condition des moments).
\item Pour tout $\nu \in K$, il existe une $1$-cha\^\i ne holomorphe
  $[S_\nu]$ de masse finie de $\Delta_n \times \CC^m \backslash \Gamma$ 
telle que $d[S_\nu]=[\Gamma,\pi,\nu]$.
\item Il existe une $p$-cha\^{\i}ne holomorphe $[T]$ de 
$\Delta_n \times \CC^m \backslash \Gamma$ dont la projection sur
$\CC^m$ est born\'ee et telle que $d[T]=[\Gamma]$.\par
\end{enumerate}
\end{prop}
\begin{proof}
L'\'equivalence entre $1.$ et $2.$ provient directement de \cite{dinh}. 
Pour toute $\phi$ fix\'e, ${\mathcal M}(z,\phi)$ est holomorphe en
$z$. L'ensemble des z\'eros de ${\mathcal M}(z,\phi)$ est donc une
hypersurface analytique $H_\phi$ de $\Delta_n$ contenant $K$. 
L'intersection $H$ de tous les $H_\phi$ est donc un ensemble analytique
contenant $K$. Comme $H$ contient $K$ il ne peut \^etre de dimension
inf\'erieure ou \'egale \`a $(n-1)$. On a donc $H=\Delta_n$.
De mani\`ere identique \`a \cite{dinh}, on montre que $\Gamma$ 
admet une solution $T$ au probl\`eme du bord dans $U$.
\end{proof}

\begin{rema}
L'hypoth\`ese $[\Gamma]\wedge \Theta=0$ pour toute $2$-forme verticale
$\Theta$ implique que pour ${\mathcal H}_{2n+1}$-presque tout $z \in
\Gamma$ tel que la multiplicit\'e de $[\Gamma]$ en $z$ soit non nulle,
l'espace tangent \`a $\Gamma$ en  $z$ ne contient pas de droite
verticale. En particulier $[\Gamma]$ ne contient pas de ``composantes verticales''.
\end{rema}

\subsection {Propri\'et\'es de convergence des suites de cha\^{\i}nes holomorphes.}

Nous nous int\'eressons ici \`a la convergence 
des suites de surfaces de Riemann de volumes uniform\'ement born\'es 
au sens suivant:

\begin{defi}
Soit $\{E_j\}$ une suite de sous-ensembles d'un espace m\'etrique
$X$. On dit que la suite $\{E_j\}$ converge au sens de Hausdorff 
(ou au sens de Bishop si les $E_j$ sont des surfaces de Riemann)
vers un ensemble $E \subset X$ (on note $E_j \rightarrow E$) si
\begin{enumerate}
\item l'ensemble $E$ co{\"\i}ncide avec l'ensemble limite de la suite
  $\{E_j\}$, i.e. est form\'e de l'ensemble des points de
la forme $lim_{j_{\nu}} x_{j_\nu}$, $x_{j_\nu} \in E_{j_\nu}$
  (en particulier, $E$ est ferm\'e dans $X$).
\item Pour tout compact $K \subset E$ et pour tout $\epsilon > 0$ il
  existe un indice $j(\epsilon,K)$ tel que $K$ appartient au
  $\epsilon$-voisinage de $E_j$ dans $X$ pour tout $j > j(\epsilon,K)$.
\end{enumerate}
\end{defi}

\begin{prop} \cite{bishop} (voir aussi \cite{chirka})
Soit $\{A_j\}$ un suite de sous-ensembles analytiques de dimension
pure $1$  d'une vari\'et\'e complexe $X$ dont le volume est
uniform\'ement born\'e sur tout compact:
$${\mathcal H}_{2}(A_j \cap K) \leq M_K < \infty \mbox{ pour tout } K \subset \subset X$$
et admettant un point d'accumulation dans $X$. 
Alors il existe une suite extraite $\{A_{\phi(j)}\}$ de la suite
$\{A_j\}$ qui converge dans $X$ 
vers un sous-ensemble analytique $A$ de dimension pure $1$.
De plus on a l'in\'egalit\'e:
$$\underline{lim}_{j \rightarrow \infty}{\mathcal H}_2(A_{\phi(j)}\cap
U) \geq
{\mathcal H}_2(A \cap U)$$
pour tout ouvert $U \subset \subset X$.
\end{prop}

Dans le cas des cha\^{\i}nes holomorphes, ce th\'eor\`eme se
g\'en\'eralise aussi:
\begin{prop}\cite{harvey3}
Soit $\{[A_j]\}$ une suite de $1$-cha\^{\i}\-nes holo\-morphes d'une
va\-ri\'et\'e complexe $X$ dont le volume est uniform\'ement
born\'e sur tout compact. Alors il existe une suite extraite  
$\{[A_{\phi(j)}]\}$ de la
suite $\{[A_j]\}$ convergeant vers une $1$-cha\^{\i}ne holomorphe $[A]$ de $X$
au sens de la topologie plate.
\end{prop}

\begin{prop}
Soit $\omega$ une vari\'et\'e complexe munie d'une distance $d$ et
$K \subset \omega$ un compact. Soit $V_\epsilon=\{ z \in \omega,
d(z,K) \leq \epsilon\}$ un $\epsilon$-voisinage de $K$. Alors
il existe une constante $C_\epsilon^K$ telle que tout ensemble analytique $A$ de
$\omega \backslash K$  irr\'eductible, de dimension $1$ et tel que
$\overline A$ contient un point $x \in K$ et un point $y \not \in V_\epsilon$
v\'erifie: $\mbox{Vol} A \geq C_\epsilon^K$.
\end{prop} 
\begin{proof}
Soit $M=\{z \in \omega, d(z,K)=\frac{\epsilon}{2}\}$. Le compact
$M$ d\'econnecte $\omega$.
Il existe donc un point $w \in M \cap A$.
Notons $B(w,\frac{\epsilon}{2})$ la boule de rayon
$\frac{\epsilon}{2}$ et de centre $w$. L'intersection $A \cap B(w,\frac{\epsilon}{2})$
est donc un sous-ensemble analytique ferm\'e et non vide de $B(w,\frac{\epsilon}{2})$.
Or pour tout compact d'une vari\'et\'e complexe et pour tout
$\epsilon > 0$,  il existe une constante $C_\epsilon$ 
minorant le volume des sous-ensembles analytiques passant par le
centre de toute boule de rayon $\epsilon$ dont l'intersection avec 
ce compact est non vide.
D'o\`u le r\'esultat.
\end{proof}

\subsection{Voisinage de Stein des surfaces de Riemann \`a bord.}

D'apr\`es un r\'esultat de Siu \cite{siu}, toute surface de Riemann
ouverte dans un espace complexe $X$ admet un voisinage de Stein.
Dans \cite{mialach}, il est montr\'e que si $S$
est une surface de Riemann \`a bord $\gamma$ connexe et lisse 
de l'espace projectif et telle que $S \cup \gamma$ ne contient aucune
surface de Riemann compacte alors $S \cup \gamma$ admet lui aussi un
voisinage de Stein. Dans \cite{dinh3}, ce r\'esultat est
\'etendu au cas des surfaces de Riemann \`a bord $C^1$ par morceaux  
et incluses dans un espace complexe $X$ quelconque. Nous remarquons
ici que la preuve donn\'ee dans \cite{dinh3} est 
encore valide pour les surfaces de Riemann dont le bord est de classe
$A_1$:
\begin{prop}
Soit $X$ une vari\'et\'e complexe. 
Soit $\gamma \subset X$ un compact de classe $A_1$
et $S$ un sous-ensemble analytique de dimension $1$ de $X \backslash
\gamma$. Supposons que $S \cup \gamma$ ne contient aucune surface de
Riemann compacte de $X$. Alors il existe un ouvert de Stein $V \subset
X$ voisinage de $S \cup \gamma$.
\end{prop}


\section{Probl\`eme du bord dans les vari\'et\'es produit.}

\subsection{Th\'eor\`eme principal.}
Soient $U$ une vari\'et\'e complexe connexe de dimension $n$,
$\omega$ une vari\'et\'e k\"ahl\'erienne disque convexe 
et $\pi:U \times \omega \rightarrow U$ la projection $(z,w)
\mapsto z$.

Dans toute la suite, nous supposerons que 
$[\Gamma]$ est un courant rectifiable, ferm\'e et maximalement
complexe de $U\times \omega$ dont le support $\Gamma$ est de classe
$A_{2n+1}$. Nous supposerons de plus, 
que pour tout compact $K \subset U$,  la projection 
de $\Gamma \cap (K \times \omega)$ sur $\omega$ est relativement
compacte dans $\omega$. 
Pour tout $z \in U$, nous noterons $[\gamma_z]$ la tranche
$[\Gamma,\pi,z]$ quand elle est bien d\'efinie. Nous noterons aussi
$\gamma_z=\Gamma \cap \{z\}\times \omega$. Le support de $[\gamma_z]$
est bien s\^ur inclus dans $\gamma_z$ mais la r\'eciproque n'est pas
toujours vraie.

\begin{defi}
Un sous-ensemble $K \subset U$ sera dit
{\it $(n-1)$-g\'en\'erique } s'il n'est pas inclus dans une r\'eunion
d\'enombrable  d'ensembles analytiques de dimension $(n-1)$
immerg\'es dans $U$.
\end{defi}

\begin{theo}
Supposons que pour tout $z \in U$ la tranche
$[\gamma_z]=[\Gamma,\pi,z]$ 
est bien d\'efinie et est un courant rectifiable ferm\'e dont le
support est de classe $A_1$.
Alors les propositions suivantes sont \'equivalentes:
\begin{enumerate}
\item Il existe un sous-ensemble $(n-1)$-g\'en\'erique $Z$ de $U$ 
tel que pour tout $z\in Z$, il existe une $1$-cha\^{\i}ne holomorphe
  $[S_z]$
de $(\{z\}\times \omega) \backslash \gamma_z$ telle que
  $d[S_z]=[\gamma_z]$.
\item Il existe un ouvert non vide $O \subset U$ tel que $[\Gamma]$ admette une
  solution au probl\`eme du bord dans $O \times \omega$ (i.e. il existe
une $(n+1)$-cha\^{\i}ne
  holomorphe $[T]$ de $(O \times \omega)\backslash \Gamma$ telle que 
$d[T]=[\Gamma]$ dans $O \times \omega$).
\item Il existe un ferm\'e $G \subset U$ de mesure de Hausdorff
  $(2n-1)$-dimensionnelle nulle
tel que $\forall z \not \in G$, il
  existe un voisinage $V_z$ de $z$ dans $U$ tel que le probl\`eme du
  bord pour $[\Gamma]$ soit r\'esoluble dans $V_z \times \omega$.
\item Il existe un ferm\'e $F$ de mesure nulle dans $U$ et une
$(n+1)$-cha\^{\i}ne holomorphe $T$ de $((U\backslash F)\times \omega)
\backslash \Gamma$ telle que $d[T]=[\Gamma]$ dans $(U\backslash
F)\times \omega$.
\end{enumerate}
\end{theo}

{$(1) \Rightarrow (2)$ : R\'eduction \`a un voisinage de Stein.} \par

Dans cette partie, les hypoth\`eses additionnelles sur les tranches
$[\gamma_z]$ donn\'ees dans l'\'enonc\'e du th\'eor\`eme 3.2 ne sont
pas n\'ecessaires et peuvent \^etre affaiblies : 

\begin{prop}
Supposons que pour toute $2$-forme lisse $\Theta$ d\'efinie sur
$\omega$, on aie $[\Gamma]\wedge \Theta=0$ et qu'il existe
un sous-ensemble $(n-1)$-g\'en\'erique $Z$ de $U$ tel que:
\begin{enumerate}
\item Pour tout $z \in Z$, la tranche $[\gamma_z]$ 
est bien d\'efinie et est un $1$-courant rectifiable dont le support est
de classe $A_1$.
\item Pour tout $z \in Z$, il existe une $1$-cha\^{\i}ne holomorphe $[S_z]$
de $(\{z\}\times \omega)\backslash \gamma_z$ telle que $d[S_z]=[\gamma_z]$.
\end{enumerate} 
Alors il existe un point $w \in Z$, un voisinage $W$ de $w$ dans
$U$ tel que le probl\`eme du bord pour $\Gamma$ admette une solution
dans $W \times \omega$ (i.e. il existe
une $(n+1)$-cha\^{\i}ne holomorphe $[T]$ de $(W\times\omega)
\backslash \Gamma$ telle que $d[T]=[\Gamma]$ dans $W \times \omega$).
\end{prop}
\begin{proof}
L'ensemble $Z$ \'etant $(n-1)$-g\'en\'erique, il existe un point $z
\in Z$ tel que l'intersection de $Z$ avec tout voisinage de $z$ est
encore $(n-1)$-g\'en\'erique.
Soit $V$ un voisinage de Stein de $\gamma_z$ obtenu gr\^ace \`a la
proposition 2.8. Soit  $V_\epsilon \subset \subset V$ 
un $\epsilon$-voisinage de
$\gamma_z$ et  $V_{\epsilon/2} \subset \subset V_{\epsilon}$ un $\epsilon/2$-voisinage de $\gamma_z$.
 Comme $\omega$ est disque convexe, il existe un compact $\widehat V_{\epsilon/2}$ tel que toute surface de
 Riemann irr\'eductible dont le bord est dans $V_{\epsilon/2}$ soit
 incluse dans  $\widehat V_{\epsilon/2}$.
Soit $U_{\epsilon/2}$ un voisinage de $z$ dans $U$ tel que pour tout $x \in
U_{\epsilon/2}$, $\gamma_x \subset V_{\epsilon/2}$.
Soit $C=C_{\epsilon/2}^{\overline V_{\epsilon/2}}>0$ la constante d\'efinie dans la proposition
2.7.
\begin{lemm}
Pour tout $x \in Z \cap U_{\epsilon/2}$, il existe une $1$-cha\^{\i}ne
holomorphe $[S_x]$ solution au probl\`eme du bord pour $[\gamma_x]$
dont le volume du support est \'egal \`a:
$$I= inf_{\{[S]; d[S]=[\gamma_x]\}} Vol S.$$
\end{lemm}
\begin{proof} 
 Soit $[S]=\sum m_i[S_i]$ une solution au probl\`eme du bord pour
 $[\gamma_x]$ o\`u les $S_i$ sont les composantes
irr\'eductibles de $S \backslash \gamma_x$.
Pour tout $\nu$, et $W_\nu$ un $\nu$ voisinage de $\gamma_x$, posons
$$[S^\nu]=\sum_{i; S_i \subset W_\nu} m_i[S_i].$$
Soit $\nu$ tel que $$Vol S^\nu < \frac{C}{4}.$$
D'apr\`es la proposition 2.7, il existe un nombre fini de composantes
irr\'eductibles \  
$\{S^{\nu}_i\}_{i=1..N}$ de $S^\nu$ non incluses dans $W_\nu$.
Pour toute composante irr\'eductible $S_i^\nu$ notons, si elle existe,
$\widetilde S_i^\nu$ une
surface de Riemann compacte irr\'eductible de $\{x\}\times \omega$ contenant
$S_i^\nu$ et $\widetilde S_i^\nu=\emptyset$ sinon.
Posons $$A=S \cup (\bigcup_{i=1..N}\widetilde S_i^\nu)$$ et 
montrons que le support de toute
solution au probl\`eme du bord $[R]$ pour $[\gamma]$ 
v\'erifiant $Vol R \leq I +\frac{C}{2}$
est inclus dans $A$. 
On a $[R]-[S]=\sum_{i=1}^{M} n_i [V_i]$ o\`u 
les $[V_i]$ sont les courants d'int\'egration sur des surfaces de
Riemann compactes irr\'eductibles $V_i$.
Si $V_i$ contient
une composante irr\'eductible $S_j^{\nu}$ de 
$S \backslash \gamma_x$ non incluse
dans $W_\nu$, on a par d\'efinition de $\widetilde S_j^{\nu}$,
$V_i \subset \widetilde S_j^{\nu} \subset A$. 
Donc, si pour tout $i \in \{1,...,M\}$, $V_i$ contient
une composante irr\'eductible $S_j^{\nu}$ de 
$S \backslash \gamma_x$ non incluse
dans $W_\nu$, on a $Supp([R]-[S]) \subset A$ et donc $R \subset A$.
Dans le cas contraire, il existe un $i_0 \in
\{1,...,M\}$ tel que $V_{i_0}$ ne contient que des composantes
irr\'eductibles de $S\backslash \gamma_x$ incluse dans $W_\nu$.
Soit $[\widetilde R]$ (resp. $[\widetilde S]$) 
la restriction de $[R]$ (resp. de $[S]$) \`a $V_{i_0}$.
On a alors par construction $$Vol \widetilde S \leq Vol S^\nu < \frac{C}{4}$$
Or $[\widetilde R]=[\widetilde S]+n_{i_0}[V_{i_0}]$, et donc
$$Vol \widetilde R \geq Vol V_{i_0} - Vol \widetilde S > \frac{3C}{4}.$$
Le courant $[R]-[\widetilde R]+[\widetilde S]$ est encore solution au
probl\`eme du bord pour $[\gamma]$ et v\'erifie
$$Vol Supp([R]-[\widetilde R]+[\widetilde S])
=Vol R - Vol \widetilde R +Vol \widetilde S < I + \frac{C}{2} -
\frac{3C}{4} +\frac{C}{4} \leq I$$
ce qui donne la contradiction recherch\'ee et prouve que $R\subset A$. 
L'ensemble $A$ ne contenant qu'un ensemble fini de surfaces de Riemann
compactes, on est r\'eduit au cas o\`u $\gamma$ est inclus dans
une surface de Riemann irr\'eductible $L$.
Si $[R_1]$ et $[R_2]$ sont deux solutions distinctes au
probl\`eme du bord pour $[\gamma]$ dans $L$, on a
$[R_1]-[R_2]=m[L]$ pour un certain $m \in \ZZ$. En particulier,
si deux solutions sont nulles sur une composante
irr\'eductible de $L \backslash \gamma$, elles sont \'egales.
Il existe donc une solution au probl\`eme du bord pour $[\gamma]$
dont le compl\'ementaire du support est de volume maximal et
donc tel que le support est de volume minimal.
\end{proof}

Pour tout $x \in Z$ choisissons $[S_x]$, 
une $1$-cha\^{\i}ne holomorphe solution au
probl\`eme du bord pour $[\gamma_x]$ minimisant le volume du support
des solutions.
Posons $$E_{k,l}=\left\{x \in Z \cap U_{\epsilon/2} , \mbox{Vol}[S_x] \leq k;
\frac {l C}{2} < \mbox{ Vol }S_x \leq  \frac{(l+1) C}{2} \right\}$$
L'ensemble $Z \cap U_{\epsilon/2}$ 
\'etant $(n-1)$-g\'en\'erique, 
il existe $k,l \in \NN$ tels que $E_{k,l}$ soit encore 
$(n-1)$-g\'en\'erique. 

\begin{lemm}
Soit $\{x_i\}$ une suite de points de $E_{k,l}$ convergeant vers
un point $x_\infty \in E_{k,l}$. Alors il existe une suite
extraite (que l'on notera encore $\{x_i\}$) tel que
\begin{enumerate}
\item La suite de $1$-courants
rectifiables $\{[\gamma_{x_i}]\}$ converge au sens des courants
vers $[\gamma_{x_\infty}]$.
\item La suite de supports $S_{x_i}$ converge au sens de Bishop
vers une surface de Riemann $S_\infty$ telle que $S_\infty \cup
\gamma_{x_\infty}$ ne contient aucune surface de Riemann compacte.
\end{enumerate}
\end{lemm}
\begin{proof}
La masse des $1$-cha\^{\i}nes holomorphes $[S_{x_i}]$ \'etant uniform\'ement
born\'ee, il existe une suite-extraite qui converge vers une
$1$-cha\^{\i}ne limite $[S]$.
La suite des bords $[\gamma_{x_i}]=d[S_{x_i}]$ converge alors
vers un courant limite rectifiable $[\gamma_\infty]=d[S]$ dont le
support est inclus dans $\gamma_{x_\infty}=\Gamma \cap \{x_\infty\}\times
\omega$.
D'apr\`es la
proposition 2.1, pour toute $(1,0)$-forme $\phi$ d\'efinie dans $V$ et
holomorphe au voisinage de $\gamma_{x_\infty}$, on a:
$$([\gamma_\infty],\phi)=\lim_{i \rightarrow \infty}([\gamma_{x_i}],\phi)= 
\lim_{i \rightarrow \infty} {\mathcal M}(x_i,\phi)=
{\mathcal M}(x_\infty, \phi)=
([\gamma_{x_\infty}],\phi)$$
En particulier, le courant $[R]=[\gamma_\infty] -[\gamma_{x_\infty}]$
est orthogonal \`a toutes les fonctions m\'eromorphes dont les p\^oles
ne rencontrent pas $\gamma_{x_\infty}$. Du fait que $\gamma_{x_\infty}$ est
de mesure de Hausdorff deux dimensionnelle nulle, 
on d\'eduit qu'il est rationnellement convexe dans $V$ et donc
que $[R]=0$ d'apr\`es le lemme suivant 
ce qui termine la preuve de la propri\'et\'e (1).\par
\begin{lemm}
Soit $[R]$  un $1$-courant de masse finie et de 
support inclus dans un compact $K \subset \CC^m$ (pour un certain $m
\in \NN)$) de mesure de Hausdorff $2$-dimensionnelle nulle et 
tel que $([R],\phi)=0$ pour toute les $(1,0)$-formes $\phi$ 
holomorphes au voisinage de $K$. Alors $[R]=0$.
\end{lemm}
\begin{proof}
Ce r\'esultat est bien connu dans $\CC$. Pour le cas g\'en\'eral,
il suffit de remarquer que si $K$ est de mesure de Hausdorff
$2$-dimensionnelle nulle, alors toute projection de $K$ le sera aussi.
Pour toute droite affine $\CC_\nu \subset \CC^m$
l'image direct de $[R]$ par la projection orthogonale sur $\CC_\nu$
est alors nulle ce qui montre que $[R]=0$.
\end{proof}

Pour tout $x \in E_{k,l}$, on a $S_x \subset \widehat V_{\epsilon/2}$.
Si pour tout $i \in \NN, S_{x_i} \subset V_{\epsilon} \subset V$
la limite au sens de Bishop est encore incluse dans $V$ qui ne
contient pas de courbes compactes et le lemme est v\'erifi\'e.
Sinon, on peut supposer que $\forall i \in \NN, S_{x_i} \not \subset V_{\epsilon}$.
Donc,  toutes les surfaces $S_{x_i}$ contiennent un point $w_i \not \in
V_\epsilon$ et $w_i \in \widehat V_{\epsilon/2}$.
Quittes \`a extraire une sous-suite de $\{[S_{x_i}]\}$, 
on peut aussi supposer que la suite des supports $\{S_{x_i}\}$ converge au sens de Bishop vers une surface de
Riemann limite $S_\infty$). 
On a alors:
$$\frac{l C}{2} \leq \mbox{ Vol } S_{x_\infty} \leq \mbox{Vol} S  \leq
\mbox{Vol} S_\infty  
\leq \underline {lim}_{j \rightarrow \infty}
\mbox{ Vol }S_{x_i} \leq  \frac{(l+1) C}{2}.$$
Car on a bien s\^ur $S \subset S_\infty$. 
Montrons que $S_\infty \cup
\gamma_{x_\infty}$ ne contient aucun sous-ensemble analytique compacte
$L$ de $\{x_\infty\} \times \omega$. En effet, si c'\'etait le
cas, on aurait $L \subset \widehat V_{\epsilon/2}$ car $L$ est
dans l'adh\'erence des $S_{x_i}$. Notons $\{L_i\}_{i \in \NN}$ les composantes
irr\'eductibles de $L \backslash \gamma_{x_\infty}$. 
On a alors la d\'ecomposition $[S]=\sum m_i [L_i] + \sum n_i [S_i] $ 
o\`u $m_i \in \ZZ$, $n_i \in \ZZ^*$ et $S_i$ sont les composantes
irr\'eductibles de $S$ non incluses dans $L$. L'ensemble analytique
$L$ n'est pas inclus dans $V$ car $V$ est une vari\'et\'e de Stein
et $L$ est compact. Il existe donc $i_0 \in \NN$ tel que $L_{i_0}$ ne soit
pas inclus dans $V$. D'apr\`es la proposition 2.6, $\mbox{Vol} L_{i_0} >
C$. La $1$-cha\^{\i}ne holomorphe
$[S^*]=[S]-m_{i_0}[L]$ est solution au probl\`eme du bord pour
$\gamma_{x_\infty}$ et a un support $S^*$ inclus dans $S_\infty
\backslash L_{i_0}$ et v\'erifie donc 
$$\mbox{Vol} S^* \leq \mbox{Vol} S_\infty-\mbox{Vol} L_{i_0} 
< \frac{(l+1)C}{2}-C \leq \frac{(l-1) C}{2} < Vol S_{x_\infty}.$$
Ce qui contredit le fait que $[S_{x_\infty}]$ minimise le support des
solutions au probl\`eme du bord pour $[\gamma_{x_\infty}]$.
\end{proof}

Soit $w$ un point de $E_{k,l}$ tel que $E_{k,l}$ soit
$(n-1)$-g\'en\'erique dans tout voisinage de $w$.
 De la m\^eme mani\`ere que dans le
lemme 3.4, il existe un sous-ensemble analytique $A$ de volume fini de
$(\{w\}\times \omega) \backslash \gamma_w$ tel que pour toute
suite $\{x_i\}$ de points de $E_{k,l}$ convergeant vers $w$,
la suite $\{[S_{x_i}]\}$ converge
vers une $1$-cha\^{\i}ne holomorphe $[S]$ solution au probl\`eme
du bord pour $[\gamma_w]$ tel que $S \subset A$.
L'ensemble $A \cup \gamma_w$ contient un nombre fini $A_1,...,A_k$ de
surfaces de Riemann irr\'eductibles et compactes. Pour tout $i=1..k$, 
l'ensemble $A_i \backslash \gamma_w$ se d\'ecompose lui aussi en
composantes irr\'eductibles $\{A_i^j\}_{j \in \NN}$.
Soit $V_{(j_1,...,j_k)}$ un voisinage de Stein de $\gamma_w \cup (A \backslash 
(\cup_{i=1..k}A_{i}^{j_i}))$. Soit $E_{k,l}^{j_1,...,j_k}=\{x \in
E_{k,l}; S_x \subset V_{(j_1,...,j_k)}\}$. 
Il existe $(j_1,...,j_k) \in
\NN^k$ tel que $E_{k,l}^{j_1,...,j_k}$ est $(n-1)$-g\'en\'erique.
En effet, dans le cas contraire,  l'ensemble union des $E_{k,l}^{j_1,...,j_k}$ serait inclus
dans une r\'eunion d\'enombrable $Y$ d'ensembles analytiques immerg\'es
de dimension $(n-1)$ de $U$ et l'ensemble $E_{k,l} \backslash Y$
serait encore $(n-1)$-g\'en\'erique. Il existerait donc une suite de
points $\{x_i\}$ de $E_{k,l}\backslash Y$ convergeant vers $w$. Quitte \`a
extraire une sous-suite, on peut supposer que la suite $\{x_i\}$
v\'erifie les propri\'et\'es $1$ et $2$ du lemme $3.5$ .
Pour $k$ assez grand, il existe alors
$(i_1,...,i_k) \in \NN^k$ tel que $x_k \in E_{k,l}^{j_1,...,j_k}$ ce
qui donne la contradiction recherch\'ee. Il existe donc un ouvert
de Stein $V_{i_1,...,i_k}$ tel que le probl\`eme du bord pour $\Gamma$
soit r\'esoluble dans $V_{i_1,...,i_k}$ 
pour un ensemble $(n-1)$-g\'en\'erique de
tranches de $\Gamma$. D'apr\`es la proposition 2.2, 
le probl\`eme du bord est r\'esoluble pour
$\Gamma$ dans $V_{i_1,...,i_k}$. 
\end{proof}
\par
\medskip


{$(2) \Rightarrow (3)$: Contr\^ole du volume jusqu'au bord}.

Supposons que le probl\`eme du bord pour $[\Gamma]$ admet une solution
$[T]$ dans $V \times \omega$ o\`u $V$ est un ouvert (assez petit) de $U$.
Alors, il existe une majoration du volume des $1$-cha\^{\i}nes holomorphes 
$[S_z]=[T,\pi,z]$  ne d\'epend que de la g\'eom\'etrie
 de $\Gamma$:
\begin{lemm} 
Soit $(x_1,...,x_{2n})$ un syst\`eme de coordonn\'ees r\'eelles de
$V$. Soit $\Pi:V \rightarrow \RR^{2n-1}$ la projection 
$\Pi(x_1,...,x_{2n})=(x_2,...,x_{2n})$. Alors pour presque tout
point $X=(x_2,...,x_{2n}) \in \RR^{2n-1}$ et pour tous points
$a,b \in \RR$ on a les in\'egalit\'es suivantes:
 $$|Vol([S_{(a,X)}])-Vol([S_{(b,X)}])|\leq M_V.Vol(\ind_{[(a,X),(b,X)]}[\Gamma,\Pi,X]) \leq
M_V.Vol([\Gamma,\Pi,X])$$
$$|Vol(S_{(a,X)})-Vol(S_{(b,X)})|\leq M_V.Vol(\Gamma\cap
[(a,X),(b,X)]\times \omega)$$
o\`u $M_V=\sup_{(V\times \omega) \cap \Gamma}|\Omega|$ pour $\Omega$ la forme k\"ahl\'erienne
associ\'ee \`a la vari\'et\'e $\omega$ et $\ind_{[(a,X),(b,X)]}$ est 
l'indicatrice de $[(a,X),(b,X)]\times \omega$ dans
$V\times \omega$.
\end{lemm}

\begin{proof}
Soit $[T]=\sum n_i [T_i]$, $n_i \in \ZZ$, on a alors  
par d\'efinition, $[\Gamma]=\sum n_i d[T_i]$.
Soit $[T^+]=\sum |n_i| [T_i]$ et $[\Gamma^+]=\sum |n_i|d[T_i]$
on a alors bien s\^ur $$Vol [\Gamma^+] \leq Vol [\Gamma].$$
D'apr\`es le th\'eor\`eme de tranchage des courants, 
pour presque tout $X\in \RR^{2n-1}$, la tranche $[\Gamma^+,\Pi,X]$
(resp. $[T^+, \Pi,X]$) est bien d\'efinie et est un $2$-courant 
(resp. $3$-courant) de masse finie.
Alors les courants $\ind_{[(a,X),(b,X)]}[\Gamma^+,\Pi,X]$ et \\
$\ind_{[(a,X),(b,X)]}[T^+,\Pi,X]$ sont bien d\'efinis car les courants $[\Gamma^+,\Pi,X]$
et $[T^+,\Pi,X]$ sont de masse localement finie.
On a alors
$$d(\ind_{[(a,X),(b,X)]}[T^+,\Pi,X])=\ind_{[(a,X),(b,X)]}d[T^+,\Pi,X]+
[T^+,\pi,(b,X)]-[T^+,\pi,(a,X)].$$
En remarquant que $d\Omega=0$ on obtient
$$0=(d(\ind_{[(a,X),(b,X)]}[T^+,\Pi,X]),\Omega)=$$
$$(\ind_{[(a,X),(b,X)]}d[T^+,\Pi,X],\Omega)+
([T^+,\pi,(b,X)],\Omega)-([T^+,\pi,(a,X)],\Omega)$$
Or par d\'efinition,   $d[T^+,\Pi,X]=[\Gamma^+,\Pi,X]$
et $([T^+,\pi,(a,X)],\Omega)$ et
$([T^+,\pi,(b,X)],\Omega)$ sont les masses des courants $[S_{(a,X)}]$
et $[S_{(b,X)}]$. On a donc
$$|Vol [S_{(a,X)}] - Vol [S_{(b,X)}]|\leq M_V.Vol (\ind_{[(a,X),(b,X)]}[\Gamma^+,\Pi,X]) \leq$$
$$M_V. Vol (\ind_{[(a,X),(b,X)]}[\Gamma,\Pi,X])$$
 ce qui montre la premi\`ere in\'egalit\'e.
La deuxi\`eme in\'egalit\'e s'obtient alors en
appliquant ce qui pr\'ec\`ede au courant $[\widetilde T]=\sum [T_i]$
dont le bord $[\widetilde \Gamma]=\sum d[T_i]$ est un courant rectifiable
\`a support dans $\Gamma$ et de
multiplicit\'e $0,\pm 1$ en presque tout point de $\Gamma$.
\end{proof}

Soit $V_{max}$ l'ensemble des points $z \in U$ pour lesquels il existe un
voisinage $V_z$ de $z$ dans $U$ tel que le probl\`eme du bord pour
$[\Gamma]$ soit r\'esoluble dans $V_z \times \omega$. Bien s\^ur
$V_{max}$ est un ouvert de $U$.

\begin{lemm}
Soit $z \in V_{max}$ et $V_z$ un voisinage de $z$ dans $U$ tel
que le probl\`eme du bord pour $[\Gamma]$ admette une solution $[T]$
dans $V_z\times \omega$.
Alors il existe un ouvert connexe maximal $U_z \subset
  V_{max}$ tel que le probl\`eme du bord pour 
$[\Gamma]$ aie une solution  $[T_{max}]$ dans $U_z \times \omega$
co\"{\i}ncidant avec $[T]$ dans $V_z \times \omega$.
De plus, $\partial U_z \backslash V_{max}$ est un ferm\'e
de mesure de Hausdorff $(2n-1)$-dimensionnelle nulle.
\end{lemm}
\begin{proof}
D'apr\`es le lemme de Zorn, 
pour montrer l'existence de $U_z$, il suffit de montrer que si 
$U_1$ et $U_2$ sont deux ouverts connexes de $U$
contenant $V_z$,  tels que $U_1 \subset U_2$ et tels que la solution
$[T_{U_1}]$ (resp. $[T_{U_2}]$) au probl\`eme du bord pour $[\Gamma]$
dans $U_1\times \omega$ (resp. $U_2 \times \omega$) co\"{\i}ncide avec
$[T]$ dans $V_z \times \omega$, alors $[T_{U_2}]$ co\"{\i}ncide avec
$[T_{U_1}]$ dans $U_1 \times \Omega$.
La $n$-cha\^{\i}ne holomorphe  $[L]=[T_{U_2}]-[T_{U_1}]$ est bien
d\'efinie dans $(U_1 \times \omega) \backslash \Gamma$. 
Par hypoth\`ese, $[L]$ est ferm\'ee et $[L]=0$ dans $V_z \times
\omega$. Le support $L$ de $[L]$ est donc un sous-ensemble analytique 
de $U_1 \times \omega$ dont l'intersection avec $V_z \times \omega$ est
vide. Soit $\pi:U\times \omega \rightarrow U$ la projection canonique
sur $U$. D'apr\`es le th\'eor\`eme de l'application propre, la
projection $\pi(L)$ de $L$ sur $U_1$ est un sous-ensemble analytique de $U_1$. 
Du fait que pour tout $z \in U$, la tranche $\{z\}\times \omega \cap
\Gamma$ est de classe $A_1$, 
on en d\'eduit que $\pi(L)$ et soit vide soit 
un sous-ensemble analytique de dimension $n$ de $U$ (i.e. $\pi(L)=U_1$
car $U_1$ est connexe). Comme $L \cap V_z\times \omega = \emptyset$, $L$
est n\'ecessairement vide dans $U_1 \times \omega$ et donc
$[L]=[T_{U_2}]-[T_{U_1}]=0$, d'o\`u l'existence de $U_z$.
Supposons que l'ensemble $G_z=\partial U_z \backslash V_{max}$ soit un ferm\'e
de mesure de Hausdorff $(2n-1)$-dimensionnelle non nulle.
Soit alors $x_0$ un point de $G_z$ tel que $G_z$ soit de mesure
non nulle dans tout voisinage $x_0$. Choisissons alors un
voisinage $V_{x_0}$ de $x_0$ et un syst\`eme de coordonn\'ees
$(x_1,...,x_{2n})$ de
$V_{x_0}$ tels que :
\begin{enumerate}
\item L'ensemble $H=\Pi(G_z \cap V_{x_0})$ est de mesure non
nulle dans $\RR^{2n-1}$ o\`u $\Pi:\RR_+\times \RR^{2n-1} \rightarrow
\RR^{2n-1}$ est la projection $\Pi(x_1,...,x_{2n})=(x_2,...,x_{2n})$.
\item L'ensemble $\{0\} \times \RR^{2n-1} \subset V_z$. 
\end{enumerate}
D'apr\`es le th\'eor\`eme de tranchage des courants, pour presque tout
$X \in \RR^{2n-1}$, le courant $[\Gamma,\Pi,X]$ est bien d\'efini
et est un $2$-courant rectifiable. Quittes \`a restreindre $H$,
on peut supposer que ceci est v\'erifi\'e pour tous les points de $H$.
Soit $$\widetilde H=\{(\lambda_X,X)\in \overline{\RR_+}\times H;
\lambda_X=\sup \{t\in \RR_+; [0,t[\times \{X\} \subset V_z\}\}$$
Soit $(\lambda_X,X) \in \widetilde H$, par d\'efinition de $H$, $\lambda_X < \infty$,
d'apr\`es le lemme 3.7, le volume des $1$-cha\^{\i}nes holomorphes
$[S_{(t,X)}]=[T,\pi,X]$ (o\`u $\pi:U\times \omega \rightarrow U$, $\pi(x,w)=x$)
pour $t < \lambda$ est uniform\'ement born\'e. D'apr\`es la
proposition 2.6, il existe alors une suite de r\'eels $\{t_i\}$ convergeant
vers $\lambda$ tel que la suite de $1$-cha\^{\i}nes holomorphes $[S_{(t_i,X)}]$
converge au sens des courants vers une $1$-cha\^{\i}ne holomorphe
$[S_{(\lambda,X)}]$ solution au probl\`eme du bord pour $[\gamma_{(\lambda,X)}]$.
On a donc $\widetilde H \subset G_z$ et est $(n-1)$-g\'en\'erique
car de mesure de Hausdorff $(2n-1)$-dimensionnelle non nulle.
La proposition 3.3 permet alors d'obtenir la contradiction recherch\'ee.
\end{proof}

 Soit $G=\partial V_{max}$, il reste \`a
montrer que $G$ est de mesure de Hausdorff $(2n-1)$-dimensionnelle
nulle. En effet, si c'est le cas on aura $G=U \backslash V_{max}$ 
car $U$ est connexe et un ensemble de mesure de Hausdorff
$(2n-1)$-dimensionnelle nulle ne peut disconnecter un ouvert de
dimension $2n$. Supposons que $G$ est de mesure
$(2n-1)$-dimensionnelle non nulle. Il existe alors un point $z \in G$ tel que 
pour tout voisinage $V_z$ de $z$ dans $U$, l'ensemble $V_z \cap G$ 
soit de mesure de Hausdorff $(2n-1)$-dimensionnelle non nulle.
Soit $V$ un voisinage de Stein de $\gamma_z=\Gamma
\cap \{z\}\times \omega$. Soit $C=C_{2\epsilon}^{\overline
V_\epsilon}$ la constante d\'efinie dans la proposition 2.7 avec
$\epsilon$ choisi de mani\`ere \`a ce que le voisinage $V_{2\epsilon}$ soit
inclus dans $V$. Soit $x_0 \in V_{max}$ tel que $\gamma_{x_0}
\subset V_{\epsilon}$. Soit $V_z$ un voisinage de $z$ dans $U$ et
$(x_1,...,x_{2n})$ un syst\`eme de coordonn\'ees de $V_z$ tel que:
\begin{enumerate}
\item $\forall x \in V_z$, $\gamma_x \subset V_\epsilon$. 
\item L'ensemble $H=\Pi(G \cap V_{z})$ est de mesure non
nulle dans $\RR^{2n-1}$ o\`u $\Pi:\RR_+\times \RR^{2n-1} \rightarrow
\RR^{2n-1}$ est la projection $\Pi(x_1,...,x_{2n})=(x_2,...,x_{2n})$.
\item L'ensemble $\{0\} \times \RR^{2n-1} \subset U_{x_0} \cap
V_z$ o\`u $U_{x_0}$ est un ouvert maximal contenant $x_0$ tel que le
probl\`eme du bord pour $[\Gamma]$ admette une solution $[T_{x_0}]$ dans
$U_{x_0} \times \omega$.
\end{enumerate}
D'apr\`es le th\'eor\`eme de tranchage des courants et le th\'eor\`eme
de Fubini, quittes \`a restreindre $H$, on peut supposer que pour tout $X
\in H$, le courant $[\Gamma,\Pi,X]$ est bien d\'efini, de masse finie
 et l'ensemble $\Gamma \cap ((\RR\times X)\times \omega)$ 
est de mesure de Hausdorff
$2$-dimensionnelle finie. L'ensemble $H$
\'etant $(n-1)$-g\'en\'erique, il existe $k,l \in \NN$ tel que
l'ensemble $$H_{l}^0=\{(0,X) \in V_z; Vol S_{(0,X)}<k
\mbox{ et } Vol [S_{(0,X)}]<l \}$$ (o\`u
$[S_{(0,X)}]=[T_{x_0},\pi,(0,X)]$) soit encore $(n-1)$-g\'en\'erique.
Pour tout $X \in \Pi(H_{l}^0)$, notons $\lambda_X$ la borne
sup\'erieure des r\'eels $t$ tels qu'il existe un voisinage $V_{(t,X)}$ de
$[0,t[\times \{X\}$, tels que le probl\`eme du bord pour $[\Gamma]$ admette une
solution $[T_{(t,X)}]$ dans $V_{(t,X)} \times \omega$ qui co\"{\i}ncide avec
$[T_{x_0}]$ au voisinage de $(0,X)\times \omega$.
Soient $X\in \Pi(H_{l}^0)$, $0<t_1<t_2<...<t_n<...<\lambda_X$ une suite de
r\'eels convergeant vers $\lambda_X$ et $[S_n]=[T_{(t_n,X)},\pi,(t_n,X)]$.
D'apr\`es le lemme 3.7, pour tout $X \in \Pi(H_{l}^0)$,
$$Vol [S_n] < l+M_V.Vol(\ind_{[(0,X),(\lambda,X)]}[\Gamma,\Pi,X])$$
et
$$Vol S_n < k+M_V.Vol(\Gamma\cap ([(0,X),(\lambda,X)]\times \omega)).$$
On en d\'eduit qu'il existe une suite extraite (que l'on notera
encore $\{[S_n]\}$) convergeant au sens des courants vers
une $1$-cha\^{\i}ne limite $[S^\infty]$ solution au probl\`eme du
bord pour $[\gamma_{(\lambda,X)}]$ et telle que la suite des
supports $S_n$ converge au sens de Bishop vers une surface de
Riemann limite $S_\infty$.
Si $S_\infty \subset V$, pour $n$ assez grand, $T_{(t_n,X)}$
est lui aussi inclus dans $V$ et donc d'apr\`es la proposition 2.2,
le probl\`eme du bord pour $[\Gamma]$ est r\'esoluble dans $V$
ce qui contredit que $z \in G$. Dans le cas contraire,
$\gamma_{(\lambda,X)}\cup S_\infty$ contient n\'ecessairement
une surface de Riemann compacte $L$. En effet,
sinon,  en consid\'erant un voisinage de Stein de
$\gamma_{(\lambda,X)}\cup S_\infty$, on montre gr\^ace \`a la
proposition 2.2 que $[\Gamma]$ admet une solution au probl\`eme
du bord dans $W$ qui co\"{\i}ncide avec $[T_{(t_n,X)}]$ pour $n$
assez grand ce qui contredit la d\'efinition de $\lambda_X$.
De m\^eme que pr\'ec\'edemment, il existe $m\in \ZZ$ tel
que la $1$-cha\^{\i}ne holomorphe $[S^\infty]-m[L]$ soit une
solution au probl\`eme du bord pour $[\gamma_{(\lambda_X,X)}]$
dont le support est inclus dans $S_\infty \backslash \widetilde
L$ o\`u $\widetilde L$ est une composante irr\'eductible de $L \backslash
\gamma_{(\lambda_X,X)}$ dont le volume est strictement sup\'erieur \`a $C$.
On a donc $$\forall X \in \Pi(H_{l}^0), Vol S_{(\lambda_X,X)} <
k+M_V.Vol(\Gamma\cap [0,\lambda_X]\times \omega)-C.$$
D'apr\`es la proposition 3.3, il existe $x_1 \in H_{l}^0$ tel que le probl\`eme
du bord pour $\Gamma$ admette une solution au voisinage
$V_{x_1}$ de $x_1$
et tel que $H_{l}^0$ soit encore $(n-1)$-g\'en\'erique dans tout voisinage de $x_1$.
On a bien s\^ur $V_{x_1} \subset V_{max}$. En appliquant le
m\^eme raisonnement que pr\'ec\'edemment \`a $V_{x_1}$
on construit un ensemble $(n-1)$-g\'en\'erique $H_{l_1}^1$
tel que pour tout $(\lambda_X,X) \in H_{l_1}^1$ on aie:
 $$Vol S_{(\lambda_X,X)} \leq k+M_V.Vol(\Gamma \cap
[(0,X),(\lambda_{X},X)]\times \omega -2C$$
et par r\'ecurrence, on construit une suite d'ensembles
$H_{l_n}^n$ tel que  pour tout $(\lambda_X,X) \in H_{l_n}^n$ on aie:
 $$Vol S_{(\lambda_X,X)} \leq k+M_V.Vol(\Gamma \cap
[(0,X),(\lambda_{X},X)]\times \omega -(n+1)C.$$
Mais alors pour $n$ assez grand, nous obtenons la contradiction recherch\'ee.

{{$(3) \Rightarrow (4)$ et $(4) \Rightarrow (1)$}. \'Evident.


\subsection{ Obstruction \`a l'existence d'une solution globale}
Dans le th\'eor\`eme 3.2, l'ensemble $F$ peut disconnecter $U$ et les solutions
dans les diff\'erentes composantes connexes de $U \backslash F$
peuvent ne pas se recoller en une solution globale
dans $U \times \omega$ comme le montre l'exemple suivant:\par
\noindent

{\bf Exemple.}
Soient $C(1,3)=\{z \in \CC, 1<|z|<3\}$, $C(2)=\{z \in \CC,
|z|=2\}$ et $\phi:\CC \rightarrow \CC^2$ d\'efinie par
$\phi(z)=(z,e^\frac{1}{z})$.
Soient $S=\phi(C(1,3))$, $\gamma= \phi(C(2))$ et
$\Gamma=\gamma \times P_1(\CC)$. Soit $\widetilde \Gamma$ une petite
d\'eformation de $\Gamma$ 
dans $S \times P_1(\CC)$ telle que
$\forall z \in S$, $\{z\}\times P_1(\CC) \cap \Gamma$ soit inclus
dans une courbe de classe $C^1$.\par

Pour $z \in C(1,3)$ assez loin de $\gamma$, $\{(z,e^\frac{1}{z})\} \times
P_1(\CC)$ ne rencontre pas $\Gamma$. Le th\'eor\`eme 3.2 s'applique donc
mais il est ici impossible de trouver une solution globale
dans $\CC \times (\CC \times P_1(\CC))$.\par

\subsection{ Probl\`eme du bord dans $\CC^n$ et $P_n(\CC)$.}

Le fait de r\'esoudre le probl\`eme du bord dans un espace produit
n'est pas limitatif. En effet, dans le cas d'un espace
localement feuillet\'e par des sous-ensembles analytiques, on pourra
se ramener localement au cas d'un produit 
et en d\'eduire la r\'esolution du probl\`eme du bord.
Par exemple, dans l'espace projectif
nous retrouvons les r\'esultats connus  sur la r\'esolution du
probl\`eme du bord.

\begin{coro}\cite{harvey,henkin,dinh}
Soit $[\Gamma]$ un courant rectifiable, ferm\'e, maximalement complexe
de $P_n(\CC) \backslash P_{(n-p+1)}(\CC)$ dont le support est de
classe $A_{2p-1}$.
Alors il existe une $p$-cha\^{\i}ne holomorphe $[T]$
de $(P_n(\CC) \backslash P_{(n-p+1)}(\CC))\backslash \Gamma$
 telle que $[\Gamma]=d[T]$ dans $P_n(\CC) \backslash P_{(n-p+1)}(\CC)$.
\end{coro}
\begin{proof}
Pour $p=2$, $\Gamma \subset P_n(\CC) \backslash P_{n-1}(\CC) \simeq \CC^n$.
Le probl\`eme du bord pour $[\Gamma]$ 
admet donc une solution (voir \cite{dinh}).
Pour $p >2$, soit $H$ le $(n-p+1)$ espace affine tel que 
$\Gamma \subset P_n(\CC)
\backslash H$. D'apr\`es la propri\'et\'e de tranchage des courants
rectifiables de support de classe $A_{2p-1}$ et le th\'eor\`eme de
Fubini,  on peut supposer sans
perte de g\'en\'eralit\'e qu'il existe 
une sous vari\'et\'e $V \subset G(n-p+2,n)$
form\'ee de $(n-p+2)$-plans contenant $H$ telle  
que pour presque tout $(n-p+2)$-plan $G_\nu$ avec $\nu \in V$,
$[\Gamma \cap G_\nu]$ est un $3$-courant rectifiable, ferm\'e et
maximalement complexe dans 
$G_\nu \backslash H \simeq \CC^{n-p+2}$ et dont 
le support est de classe $A_3$.
Le probl\`eme du bord pour le courant d'intersection $[\Gamma \cap G_\nu]$
admet donc une solution dans $G_\nu \backslash H$.
Soit $\Psi: V \times \CC^{n-p+2}\rightarrow P_n(\CC)\backslash H$ 
une application holomorphe injective identifiant
$\{\nu\}\times \CC^{n-p+2}$ avec $G_\nu\backslash H$
(on a suppos\'e ici que $V$ est assez petit pour avoir l'\'existence
de $\Psi$).
Soit $[\widetilde \Gamma]=\Psi_*([\Gamma])$ l'image inverse de
$[\Gamma]$ par $\Psi$. D'apr\`es la proposition 2.2, 
le probl\`eme du bord est 
r\'esoluble pour $[\widetilde \Gamma]$ dans $V \times \CC^{n-p+2}$. 
Et donc il existe un $(n-p+2)$-plan $P$ tel que le probl\`eme du
bord pour $[\Gamma]$ admette une solution $[T]$ au voisinage de $P$.
Soit $P_\epsilon$ un $\epsilon$-voisinage de $P$ dans $P_n(\CC)$.
Pour $\epsilon$ assez petit la restriction de $[T]$ \`a $P_\epsilon$ 
est une $p$-cha\^{\i}ne holomorphe $[T_\epsilon]$ dont le bord est un
courant rectifiable maximalement complexe dont le support est de
classe $A_{2p-1}$ et qui co\"{\i}ncide avec $[\Gamma]$ au voisinage de
$P$. La r\'esolution du probl\`eme du bord pour $[\Gamma]$ et alors
\'equivalente \`a la r\'esolution du probl\`eme du bord pour
$[\Gamma]-d[T_\epsilon]$ dont le support est dans $P_n(\CC) \backslash
P$.
On est donc ramen\'e au cas o\`u
$\Gamma \subset P_n(\CC) \backslash P_{n-p+2}(\CC)$. Si $p=3$, le
r\'esultat est prouv\'e. Sinon, en faisant \`a nouveau la m\^eme
manipulation, on est ramen\'e au cas
$\Gamma \subset P_n(\CC) \backslash P_{n-p+3}(\CC)$. Puis par r\'ecurrence, on se ram\`ene au
cas $\Gamma \subset P_n(\CC) \backslash P_{n-1}(\CC) \simeq \CC^n$ ce qui montre le r\'esultat.
\end{proof}

\begin{defi}[\cite{dinh2}]
Soit $U$ un ouvert de la grassmanienne $G(n-p+1,n)$, un
sous-ensemble $K \subset U$ est dit {\it projectivement $k$-g\'en\'erique}
si pour tout ensemble non $k$-g\'en\'erique $S \subset
P_n(\CC)$,  l'ensemble $K \backslash G_S$ est $k$-g\'en\'erique o\`u
$$G_S=\{\nu \in G(n-p+1,n); P_{\nu}^{n-p+1}(\CC) \cap S \neq \emptyset\}.$$
\end{defi}
Une propri\'et\'e principale des ensembles projectivement $k$-g\'en\'eriques est que
si $V\subset U$ est une sous-vari\'et\'e complexe de dimension $(p-1)$
qui est $(p-2)$-g\'en\'erique alors il existe $\nu \in V$ tel que
$P_{n-p+1}^{\nu}\cap S =\emptyset$, o\`u
$$S=\{z \in P_n(\CC); dim\{\nu \in V; z \in P_{n-p+1}^{\nu}\}\geq 1\}.$$
\begin{coro} [\cite{henkin,dinh2}]
Soit $X$ un ouvert $(n-p+1)$-concave de $P_n(\CC)$ avec $X^*=\{ \nu
\in G(n-p+1,n); P^\nu_{n-p+1} \subset X\}$ connexe. Soit $[\Gamma]$
une combinaison lin\'eaire, localement finie, \`a coefficients entiers
de courants d'int\'egration sur des sous-vari\'et\'es r\'eelles
$\Gamma_i$, orient\'ees, de classe $C^1$, de dimension $(2p-1)$ et
maximalement complexe dans $X$. Soit $K$ un sous-ensemble
projectivement $(p-2)$-g\'en\'erique de $X^*$. Supposons que pour tout
$\nu \in K$:
\begin{enumerate}
\item $P_{n-p+1}^{\nu}$ intersecte $\Gamma_i$ transversalement en tout
  point.
\item Il existe une $1$-cha\^{\i}ne holomorphe $[S_\nu]$ de masse
  finie de $P_{n-p+1}^{\nu} \backslash \Gamma$ telle que
  $d[S_\nu]=[\gamma_\nu]$ au sens des courants dans
  $P_{n-p+1}^{\nu}$ (o\`u $[\gamma_\nu]=[\Gamma\cap P_{n-p+1}^{\nu}(\CC)]$ est le
  courant d'intersection de $P_{n-p+1}(\CC)$ avec $[\Gamma])$.
\end{enumerate}
Alors il existe une $p$-cha\^{\i}ne holomorphe  $[T]$ de $X \backslash
\Gamma$, de masse localement finie, v\'erifiant 
$d[T]=[\Gamma]$ dans $X$.

\end{coro}
\begin{proof}
D'apr\`es la d\'emonstration de la proposition 3.3, 
il existe $\nu_0 \in K$, un voisinage de Stein $V_0$ de 
$S_{\nu_0} \cup \gamma_{\nu_0}$ 
(o\`u $\gamma_{\nu_0}=\Gamma \cap P_{n-p+1}^{\nu_0}(\CC)$)
et un sous-ensemble d\'enombrablement $(n-1)$-g\'en\'erique $K_0$
de $K$ tel que pour tout $\nu \in K_0$, $S_{\nu} \subset V_0$.
De la proposition 2.2, on d\'eduit alors qu'il existe une
sous-vari\'et\'e
analytique $\widetilde K \subset X^*$ de dimension $(p-1)$ 
et projectivement $(p-2)$-g\'en\'erique tel que
pour tout $\nu \in \widetilde K$, le probl\`eme du bord pour $[\gamma_\nu]$
soit r\'esoluble dans $V_0$. 
Soit $\Psi: \widetilde K \times P_{n-p+1}(\CC) \rightarrow P_n(\CC)$
la projection canonique qui identifie
$\{\nu\}\times P_{n-p+1}(\CC)$ avec $P_{n-p+1}^\nu(\CC)$.
Alors il existe $\nu \in \widetilde K$ tel que pour tout $z \in
P_{n-p+1}^{\nu}, \Psi^{-1}(\{z\})$ est fini. 
Il existe donc un voisinage $V$ de $\nu \in \widetilde K$ tel que la restriction
de $\Psi$ \`a $V \times P_{n-p+1}(\CC)$ soit propre et d'ordre fini. 
L'image inverse $[\widetilde \Gamma]=\Psi_*([\Gamma])$ est donc bien d\'efinie et est
encore un courant rectifiable, ferm\'e,  maximalement complexe
dont le support est de classe $A_{2p-1}$. D'apr\`es la proposition
3.3, il existe un ouvert $U$ de $V$ tel que le probl\`eme du bord
pour $[\widetilde \Gamma]$ admette
une solution $[\widetilde T]$ dans $U \times P_{n-p+1}(\CC)$.
Mais alors, d'apr\`es le lemme pr\'ec\'edent, l'image directe 
de $[\widetilde T]$ donne une solution $[T]=\Psi^*([\widetilde T])$ 
au probl\`eme du bord pour $[\Gamma]$ dans un voisinage
d'un $(n-p+1)$-plan $P$ de $P_n(\CC)$. Soit $P_\epsilon$ un
$\epsilon$-voisinage de $P$ dans $P_n(\CC)$. Pour $\epsilon$ assez
petit, la restriction $[T_\epsilon]$ de $[T]$ \`a $P_\epsilon$ est une
$p$-cha\^{\i}ne holomorphe dont le bord est un courant rectifiable
maximalement complexe dont le support est de classe $A_{2p-1}$ et qui
co\"{\i}ncide avec $[\Gamma]$ au voisinage de $P$. Le courant
$[\Gamma]-d[T_\epsilon]$ a donc un support dans $P_n(\CC) \backslash
P$ et admet une solution $[R]$ au probl\`eme du bord d'apr\`es
le corollaire pr\'ec\'edent. La $p$-cha\^{\i}ne holomorphe
$[R]+[T_\epsilon]$ nous donne alors une solution au probl\`eme du bord
pour $[\Gamma]$.
\end{proof}
\par

\subsection{Th\'eor\`eme de Hartogs-Levi g\'en\'eralis\'e.}

On note $\Delta$ le disque unit\'e de $\CC$, $C(r)=\{z\in \CC,
|z|=r\}$ et $C(r_1,r_2)=\{z \in \CC, r_1 <|z|<r_2\}$. 

\begin{coro}(Th\'eor\`eme de Hartogs-Levi g\'en\'eralis\'e)\par\noindent
Soit $X$ une vari\'et\'e K\"ahl\'erienne disque convexe.
Soit $f$ une application m\'eromorphe, \`a valeurs dans $X$
et d\'efinie sur $C(1-\epsilon,1)\times \Delta$. Soit $\{l_\nu\}_{\nu \in
  V}$ une famille non d\'enombrable de droites complexes.
On suppose de plus que pour tout $\nu \in V$, 
$l_\nu \cap C(1-\epsilon,1)\times \Delta
\neq \emptyset$, $l_\nu \cap C(1-\epsilon,1)\times b\Delta= \emptyset$ et
pour tous $\nu_1$, $\nu_2$, $\nu_3 \in V$ deux \`a deux diff\'erents,
 $l_{\nu_1}\cap l_{\nu_2}\cap l_{\nu_3}\cap \Delta^2= \emptyset$.
Supposons que $f$ se prolonge m\'eromorphiquement \`a 
$l_\nu \cap \Delta^2$ pour tout $\nu \in V$. Alors $f$ se prolonge
m\'eromorphiquement \`a $\Delta^2$.
\end{coro}
\begin{proof}
L'application $f$ \'etant m\'eromorphe, l'ensemble des points
d'in\-d\'e\-ter\-mina\-tion de $f$ est discret dans
$C(1-\epsilon,1)\times \Delta$. Soir $r \in ]1-\epsilon,1[$, tel que 
$M_r=C(r)\times\Delta$ ne rencontre 
aucun point d'ind\'etermination de $f$. 
La restriction de $f$ \`a $M_r$ est donc lisse et
pour tout $\nu \in V \subset G(2,3)$, la droite $l_\nu$ est transverse \`a
$M_r$.   Soit $\Gamma_f=\{(w,c)\in \Delta^2\times X, w\in M_r,
c=f(w)\}$ le graphe de $f$. 
De la m\^eme mani\`ere que dans la proposition 3.3,
il existe une surface de Riemann $V \subset G(3,2)$,
un point $\nu_0 \in V$, un
voisinage de Stein $W$ de $S_{\nu_0} \cup \gamma_{\nu_0}$
(o\`u $S_{\nu_0}$ est le graphe de l'extension m\'eromorphe de $f$ \`a $l_{\nu_0}$
et $\gamma_{\nu_0}=\Gamma_f \cap (l_{\nu_0} \times X)$),
tels que l'ensemble des points $\nu \in V$ tel que $S_{\nu} \subset W$ 
est $0$-g\'en\'erique. Soit $\phi: W \rightarrow \CC^m$ un plongement
de $W$ dans l'espace affine. D'apr\`es l'hypoth\`ese faite
sur les droite $l_{\nu}$ et le th\'eor\`eme de
 Levi, $\phi \circ f$
admet une extension holomorphe au voisine de $\Delta^2 \cap
l_{\nu_0}$. Et donc $f=\phi^{-1}\circ \phi \circ f$ admet une
extension holomorphe au voisinage de   $\Delta^2 \cap
l_{\nu_0}$. L'extension m\'eromorphe \`a tout $\Delta^2$ se fait alors soit en
appliquant \`a nouveau le th\'eor\`eme 3.2
soit gr\^ace au r\'esultat de \cite{Ivashkovich}.
\end{proof}
\par

\subsection{G\'en\'eralisations du th\'eor\`eme de Hartogs-Bochner.}

Soit  $X$ une vari\'et\'e  k\"ahl\'erienne disque convexe.
Soit $\Gamma$ une sous-vari\'et\'e, de classe $C^1$, orient\'ee, ferm\'ee et
maximalement complexe de $\CC^n$. Supposons qu'il existe un sous-ensemble
analytique et irr\'eductible $A$ de $\CC^n\backslash \Gamma$ tel
que $d[A]=[\Gamma]$ o\`u $[A]$ et $[\Gamma]$ sont les courants
d'int\'egration sur $A$ et $\Gamma$. 
Soit $\pi:\CC^n \rightarrow \CC^{p-1}$ la projection $\pi(z_1,...,z_n)=(z_1,...,z_{p-1})$.
Nous supposerons que 
pour tout $z \in \CC^{p-1}$, $\Gamma \cap \{z\}\times \CC^{n-p+1}$ est
un compact de classe $A_1$.

\begin{coro}
Toute application Lipschitzienne CR $f: \Gamma \rightarrow X$ admet une extension
m\'eromorphe \`a $A$.
\end{coro}
\begin{proof}
Soit $$\Gamma_f=\{(x,(y,w)) \in \CC^{p-1}\times(\CC^{n-p+1} \times X),
w=f(x,y)\}.$$
La vari\'et\'e $\Gamma$ \'etant compacte, il existe $R>0$ tel que si
$|x|>R$, on aie $\pi^{-1}(x) \cap \Gamma = \emptyset$. Le probl\`eme du
bord pour $[\Gamma_f]$ est donc r\'esoluble dans l'ouvert
$\{|x|>R\}\times (\CC^{n-p+1} \times X)$. La vari\'et\'e
$\omega=\CC^{n-p+1} \times X$ est bien s\^ur disque convexe, le
th\'eor\`eme 3.2 s'applique donc.
Il reste \`a v\'erifier que  l'extension ainsi obtenue donne une
solution globale au probl\`eme du bord pour $[\Gamma_f]$.
Soit $U \subset \CC^{p-1}$ un ouvert connexe maximal contenant l'ouvert
$\{|x|>R\}$ tel que le probl\`eme du
bord pour $[\Gamma_f]$ admette une solution  $[T_U]$ dans 
$U \times (\CC^{n-p+1} \times
X)$ qui est nulle pour $\{|x|>R\}$.
Par unicit\'e du probl\`eme du bord dans $\CC^n$, il est ici
impossible que $[\Gamma_f]$ admette deux solutions distinctes au
probl\`eme du bord dans $V \times \omega$ pour tout ouvert $V$ de
$\CC^{p-1}$. D'apr\`es le lemme 3.8, 
le bord $G=\partial U=\CC^{p-1}\backslash U$ 
est donc de mesure de Hausdorff $(2p-3)$-dimensionnelle nulle. 
D'apr\`es le lemme 3.7, la solution $[T]$ au
probl\`eme  du bord pour  $[\Gamma_f]$ est de volume born\'e au
voisinage de $G$, elle admet donc une extension simple sur
$\CC^{p-1}\times \omega \backslash \Gamma_f$ qui reste ferm\'ee au
sens des courants et est donc solution au probl\`eme du bord pour
$[\Gamma_f]$ dans $\CC^{p-1} \times \omega$. La projection
de cette solution sur $\CC^n$ donne alors l'extension m\'eromorphe de
$f$ sur $A$.
\end{proof}

Dans le cas o\`u $\Gamma$ est de classe $C^2$, $\Gamma$ sera une vari\'et\'e CR
globalement minimale. Les applications CR continues sur $\Gamma$
admettent une extension holomorphe sur une ``extension analytique
\`a un cot\'e'' de $\Gamma$. Quitte \`a d\'eformer $\Gamma$ 
dans cette extension,
on peut alors supposer sans perte de g\'en\'eralit\'e les hypoth\`eses
du corollaire pr\'ec\'edent v\'erifi\'ees et nous obtenons:

\begin{coro} 
Soit $\Gamma$  une sous-vari\'et\'e de classe $C^2$, compacte, connexe,  orient\'ee 
et maximalement complexe de $\CC^n$,  bord (au sens des courants) 
d'un sous-ensemble analytique born\'e irr\'eductible $A$ de 
$\CC^n \backslash \Gamma$. Alors
toute application CR continue $f: \Gamma \rightarrow X$ admet une extension
m\'eromorphe \`a $A$.
\end{coro}

\begin{rema}
Ce dernier r\'esultat peut aussi \^etre obtenu comme corrolaire du
th\'eor\`eme d'extension du type Hartogs de \cite{Ivashkovich}. 
En effet, ce dernier permet d'\'etendre $f$ m\'eromorphiquement \`a
$Reg(A)$, l'extension \`a tout $A$ peut alors par exemple \^etre obtenue 
par un  th\'eor\`eme d'extension du type Thullen (voir \cite{siu3}). 
\end{rema}

\begin{coro}
Supposons que $X$ est une vari\'et\'e K\"ahl\'erienne 
de dimension $2$ ne contenant aucune surface de Riemann compacte.
Soit $M$  une hypersurface
r\'eelle de $X$ la s\'eparant  en deux composantes connexes
$\Omega_1$ et $\Omega_2$. Alors  toute fonction holomorphe $f$ au
voisinage de $M$ admet une extension holomorphe sur $\Omega_1$ 
ou sur $\Omega_2$.
\end{coro}
\begin{proof}
De m\^eme que pr\'ec\'edemment, en consid\'erant
$[\Gamma_f]$, 
le graphe de la restriction de $f$ sur $M$ (ou une d\'eformation de
$f$), on montre qu'il existe un ouvert maximal $U_{max}$
tel que le probl\`eme du bord pour $[\Gamma_f]$  admette une solution
$[T]$ dans $U_{max}\times X$. 
Montrons par l'absurde que $U_{max}=\CC$. Supposons le contraire.
D'apr\`es le lemme 3.8, si 
$G=\CC \backslash U_{max}$ est $0$-g\'en\'erique, il existe
un point $z \in \CC$ tel que $[\Gamma_f]$
admette deux solutions $[T_1]$ et $[T_2]$ disctintes au probl\`eme du bord
dans un  voisinage $W$ de $\{z\}\times X$.
Mais alors $[R]=[T_1]-[T_2]$ est une $2$-cha\^{\i}ne holomorphe
ferm\'e et non nulle dans $W$. En particulier, pour presque
tout $x$ assez proche de $z$, $[R,\Pi,x]$ est une $1$-cha\^{\i}ne
holomorphe ferm\'ee non nulle de $\{x\}\times X$. Pour un tel $x$,
le support de $[R,\Pi,x]$ 
est donc une r\'eunion de surfaces de Riemann compactes de $\{x\}\times
X$ ce qui contredit le fait que $X$ ne contient pas de surfaces de
Riemann compactes. 
Dans le cas o\`u le ferm\'e $G$ est $0$-g\'en\'erique
(i.e. d\'enombrable), il admet des points isol\'es.
Mais alors au voisinage de ces points, 
$[T]$ admet une extension simple ferm\'ee qui reste solution
au probl\`eme du bord pour $[\Gamma_f]$, 
ce qui contredit la maximalit\'e de $U_{max}$ et termine la preuve
du corollaire.
\end{proof}

\subsection{Plongement des structures CR}
Le but de ce paragraphe est de donner une caract\'erisation
des structures CR strictement pseudoconvexes admettant une solution au probl\`eme
du bord dans une vari\'et\'e $X$ donn\'ee.
Dans le cas o\`u $X$ est de dimension $2$ (et donc $M$ une
hypersurface r\'eelle de $X$), une caract\'erisation
de nature topologique est donn\'ee dans \cite{kato}.
Dans le cas o\`u $X$ est K\"ahl\'erienne disque convexe
de dimension quelconque, la caract\'erisation suivante est valide:

\begin{prop}
Soit $M$ une sous-vari\'et\'e orient\'ee, compacte, de classe $C^2$ et
maximalement complexe de $X$ v\'erifiant l'une des trois propri\'et\'es
 suivantes:
\begin{enumerate}
\item  $M$ est plongeable dans l'espace affine et de dimension
  sup\{erieure ou \'egale \`a $3$.
 $\CC^n$.
\item $M$ est strictement pseudoconvexe et de dimension $5$.
\item $M$ est strictement pseudoconvexe, de dimension $3$ et  est le 
bord d'une vari\'et\'e complexe abstraite.
\end{enumerate}
Alors $M$ admet une solution au probl\`eme du bord.\par
R\'eciproquement, si $M$ est strictement pseudoconvexe et admet
une solution au probl\`eme du bord alors $M$ est plongeable dans
l'espace affine.
\end{prop}
\begin{proof}
D'apr\`es \cite{kohn,boutet}, la deuxi\`eme propri\'et\'e implique
automatiquement la premi\`ere. D'apr\`es \cite{kohn,henkin2}, 
la troisi\`eme propri\'et\'e implique elle aussi la 
premi\`ere. Supposons donc qu'il existe un plongement CR
$\phi: M \rightarrow \CC^m$. Soit $\widetilde{M} \subset \CC^m$, l'image 
de $M$ par $\phi$. D'apr\`es \cite{harvey}, $\widetilde M$ admet une
solution $A$ au probl\`eme du bord. L'application $\Psi=\phi^{-1}: \widetilde M
\rightarrow X$ est une application CR. D'apr\`es \cite{Joricke},
$M$ est globalement minimal (i.e. $M$ est constitu\'e d'une seule
orbite CR). La propagation de l'extension le long 
de l'orbite CR (voir \cite{Joricke,porten})
montre alors que $\Psi$ s'\'etend holomorphiquement sur une extension
analytique \`a un cot\'e de $M$ (i.e. un ensemble qui est 
une extension analytique au voisinage de chaque point de $M$, le
cot\'e peut changer). En d\'eformant $M$ dans cette extension
analytique, on peut supposer sans perte de g\'en\'eralit\'e que
$\Psi$ s'\'etend holomorphiquement sur un ensemble analytique au
voisinage de $M$. 
D'apr\`es le corollaire 3.14,
$\Psi$ s'\'etend m\'eromorphiquement sur $A$. L'image directe de $[A]$
par l'extension de $\Psi$
nous donne un sous-ensemble analytique de $X$ solution du
probl\`eme du bord pour $[M]$.
R\'eciproquement, supposons que $M$ admet une solution au probl\`eme
du bord dans $X$, d'apr\`es les arguments de \cite{henkin2},
$M$ est plongeable dans l'espace affine.
\end{proof}\par

Dans le cas o\`u $M$ est de dimension 3, 
la question se pose alors de savoir s'il est possible de s'affranchir 
de l'hypoth\`ese de plongeabilit\'e de $M$ dans  l'espace affine. 
En effet,  on ne conna\^{\i}t pas d'exemple de
structure CR plongeable dans $X$ mais non plongeable dans $\CC^n$
(voir \cite{henkin}). Le th\'eor\`eme 3.2, permet alors de montrer
que de telles vari\'et\'es (si elles existent) n'admettent aucune
fonction CR non constante:

\begin{coro}
Soit $M$ une sous-vari\'et\'e compacte et strictement pseudoconvexe 
de dimension $3$ de $X$. Alors l'une des
deux propri\'et\'es suivantes est v\'erifi\'ee:
\begin{enumerate}
\item Les fonctions CR sur $M$ sont constantes.
\item $M$ admet une solution au probl\`eme du bord dans $X$
 et $M$ est plongeable dans l'espace affine.
\end{enumerate}
\end{coro}

\begin{proof}
Supposons 
qu'il existe une fonction CR $f:M \rightarrow \CC$ non
constante. Soit 
$$\Gamma_f=\{(f(x),x) \in \CC \times X; x \in M\}$$
le graphe de $f$.
\begin{lemm}
On peut supposer sans perte de g\'en\'eralit\'e
que $f$ est de classe $C^\infty$ et
que pour tout $z\in \CC$, 
$\gamma_z=\Gamma_f \cap \{z\}\times X$ 
est inclus dans une r\'eunion finie de courbes lisses par morceaux.
\end{lemm}
\begin{proof}
Le th\'eor\`eme de Lewy permet d'\'etendre $f$
holomorphiquement sur une extension analytique $U$ du cot\'e pseudoconvexe
de $M$ (i.e. il existe un voisinage $V$ de $M$
tel que $U \cap V$ soit un sous-ensemble analytique,
irr\'eductible de dimension $2$ de $V \backslash M$ 
tel que $d[U]=[M]$ au sens des courants).   
L'extension \'etant holomorphe 
non constante, ses lignes de niveau sont des sous-ensembles
analytiques de dimension 1.  Il existe donc une d\'eformation $\widetilde M$ 
de $M$ dans $A$  telle que les lignes de niveau de 
la restriction de l'extension de $f$ sur $M$ v\'erifient
les propri\'et\'es du lemme.
La r\'esolution du probl\`eme du bord pour $M$ et 
pour $\widetilde M$ \'etant \'equivalentes, nous obtenons la
r\'eduction recherch\'ee.
\end{proof}
La fonction $f$ \'etant continue, sont module admet un maximum 
$R$ sur $M$. Le probl\`eme du bord est donc r\'esoluble
pour $[\Gamma_f]$ dans l'ouvert $\{|z|>R\}\times X$ (prendre
comme solution $[T]=0$). D'apr\`es le lemme 3.8, il existe 
un ouvert connexe maximal $U_{max}$ contenant $\{|z|>R\}$ tel
que le probl\`eme du bord pour $[\Gamma_f]$ admette une solution
dans $U_{max}\times X$ nulle pour $|z|>R$.  
Soit $[T]=\sum_{i \in I} n_i [T_i]$ la d\'ecomposition de $[T]$ 
en composantes irr\'eductibles. 
\begin{lemm}
La $2$-cha\^{\i}ne holomorphe $[T]$ est positive 
(i.e. pour tout $i \in I$, $n_i \geq 1$).
\end{lemm}
\begin{proof}
Soit $\Pi:\CC \times X \rightarrow \CC$ la projection canonique sur
$\CC$. pour tout $i \in I$, $\pi(T_i)$ est un ouvert de $\CC$.
Soit $V=U_{max}\backslash \cup_{i \in I; n_i < 0}\Pi(T_i)$.
L'ensemble $G=\partial V \cap U_{max}$ est de mesure de Hausdorff
$(2n-1)$-dimensionnelle non nulle car il s\'epare deux ouverts de
$U_{max}$. Il existe donc un point $z \in G$, tel que $\{z\}\times X$
soit transverse \`a $T$. En particulier, $[S_z]=[T,\Pi,z]$ est une
$1$-cha\^{\i}ne holomorphe positive.
Par d\'efinition de $V$, il existe $i \in I$
tel que $n_i < 0$ et $\{z\}\times X$ soit tangent \`a $\overline
T_i$. Il existe donc un point $x \in \gamma_z$ tel que $\{z\}\times X$
soit tangent en $x$ \`a $\overline T_i$ et tel qu'il existe
un voisinage $V_x$ de $x$ tel que $d[T_i]=[\Gamma_f]$ dans $V_x$.
Comme $\{z\}\times X$ est tangent \`a $\Gamma_f$ en $x$
et que $\Gamma_f$ est strictement pseudoconvexe, n\'ecessairement
$\{z\}\times X$ est langeant du c\^ot\'e concave de $\Gamma_f$
et donc $T_i$ est du c\^ot\'e convexe de $\Gamma_f$ (i.e. $d[T_i]$ 
est le courant d'int\'egration sur $[\Gamma_f]$ avec multiplicit\'e
1). Soit $n_j \geq 0$ la multiplicit\'e de la composante irr\'eductible $T_j$ 
(si elle existe, sinon poser $n_j=0$) de $[T]$ dans $U_{max}\times X$ 
v\'erifiant $d[T_j]=-[M]$ au voisinage de $x$.
Comme $[\Gamma_f]$ est le courant d'int\'egration sur la vari\'et\'e
$\Gamma_f$, il est de multiplicit\'e $1$ en tout point. 
On a donc $n_i=n_j+1>0$  ce qui donne la contradiction
recherch\'ee. 
\end{proof}
\begin{lemm}
Soit $[\Gamma]$ un $(2p-1)$-courant rectifiable ferm\'e, maximalement complexe et dont le support est de classe $A_{2p-1}$ d'une vari\'et\'e
 complexe $Y$.
Soient $[T_1]$ et $[T_2]$
deux $p$-cha\^{\i}nes holomorphes de $Y\backslash \Gamma$,
positives et solutions au probl\`eme 
du bord pour $[\Gamma]$ dans $Y$.
Supposons que $T_1 \cup \Gamma$ et $T_2 \cup \Gamma$ ne contiennent
aucun sous-ensemble analytique de dimension $p$ de $Y$. 
Alors $[T_1]=[T_2]$.
\end{lemm}
\begin{proof}
En effet, la $n$-cha\^{\i}ne holomorphe
$[T_1]-[T_2]$  est ferm\'ee dans $Y$. On a donc 
$[T_1]-[T_2]=\sum_{i \in I} n_i [L_i]$ o\`u $n_i \in \ZZ^*$ et $L_i$ sont des sous-ensemble analytiques de $Y$ de dimension pure $p$. 
Supposons que $I$ n'est pas vide. Soit  donc 
$i \in I$ et notons $[L_i^1]$ et $[L_i^2]$
les restrictions de $[T_1]$ et $[T_2]$ \`a
$L_i$, on a donc $$[L_i^1]-[L_i^2]=n_i[L_i].$$
Quittes \`a intervertir le r\^ole de $[L_1]$ et $[L_2]$, on peut
toujours supposer que $n_i > 0$.
On a alors $L_i  \not \subset L_i^1 \cup \Gamma$ et 
$L \not \subset L_i^2 \cup \Gamma$ car $T_1 \cup \Gamma$ 
et $T_2\cup \Gamma$ 
ne contiennent pas de sous-ensembles analytiques de dimension $p$ 
de $Y$. 
Soit alors $R$ une composante irr\'eductible de $L_i
\backslash \Gamma$ telle que $[R]$ soit de multiplicit\'e $0$ dans
$[L_1]$. Par hypoth\`ese, on a $$k \geq 0$$ o\`u 
$k$ est la multiplicit\'e de $[R]$ dans $[L_2]$.
On a alors $$0-k=n_i$$ et donc $n_i$ est n\'egatif ce qui donne la
contradiction recherch\'ee. 
\end{proof}

\begin{lemm}
Soit $V$ l'ensemble des points $z \in U_{max}$ tels que
 le compact $S_z \cup \gamma_z$ contient une surface de Riemann compacte.
Alors $V$ n'est pas $0$-g\'en\'erique.
\end{lemm}
\begin{proof}
Rappelons que $S_z$ est le support de la $1$-cha\^{\i}ne
holomorphe $[S_z]=[T,\pi,z]$ o\`u $\pi:\CC \times X \rightarrow \CC$ est la projection canonique sur la base $\CC$. 
La $2$-cha\^{\i}ne holomorphe $[T]$ \'etant positive,
pour tout $z \in U_{max}$, $[S_z]$ 
est une $1$-cha\^{\i}ne holomorphe positive. 
Pour tout $z \in V$, notons $\{L_k\}_{k \in K}$ l'ensemble 
des surfaces de Riemann compactes incluses dans
$S_z \cup \gamma_z$ et $\{L_k^i\}_{i \in I}$ 
les composantes irr\'eductibles de 
$L_k \backslash \gamma_z$. 
Soit $m_k^i$ la multiplicit\'e de $[L_k^i]$ dans
$[S_z]$. Notons $m_k=\inf_{i\in I} m_k^i$ 
(nous rappelons que $m_k^i \geq 0$)
et posons $[\widetilde S_z]=[S_z]-\sum m_k [L_k]$. Par construction,
$[\widetilde S_z]$ est
la $1$-cha\^{\i}ne holomorphe positive
solution au probl\`eme du bord pour $[\gamma_z]$ telle
que $\mbox{supp }([\widetilde S_z]) \cup \gamma_z$ 
ne contient pas de surface de Riemann compacte
(l'unicit\'e d\'ecoule du lemme pr\'ec\'edent). 
Supposons que $V$ est $0$-g\'en\'erique (ou de mani\`ere \'equivalente
que $\partial V$ est $0$-g\'e\-n\'e\-ri\-que).
De mani\`ere similaire \`a la proposition 3.3,
on consid\'erant les ensembles
$$E_{k,l}=\left\{ x \in \partial V \cap U_{\epsilon/2}; 
Vol[\widetilde S_x] \leq k;
\frac{lC}{2} < Vol supp([\widetilde S_x]) \leq \frac{(l+1)C}{2} \right \}$$  o\`u $U_{\epsilon/2}$ et $C$ sont 
d\'efinis de la m\^eme mani\`ere
que dans la d\'emonstration de la proposition 3.3, on montre
qu'il existe un point $z_0 \in \partial V$, 
un voisinage $V_{z_0}$ de $z_0$ 
et un ouvert de Stein $W \subset \CC \times X$ tel que le probl\`eme
du bord pour $[\Gamma_f]$ admette une solution $[\widetilde T]$
au probl\`eme du bord
dans $W \cap (V_{z_0}\times X)$. Par construction de cette solution,
on remarque qu'il existe un sous-ensemble $0$-g\'en\'erique $\widetilde V \subset \partial V$ tel que pour tout $z \in \widetilde V$, 
$[\widetilde S_z]=[\widetilde T, \pi,z]$. 
Soit alors $z \in \widetilde V$ tel que $S_z$ soit transverse 
\`a $\widetilde T$, le fait que 
pour tout $z\in V$, $[\widetilde S_z]$ est une 
$1$-cha\^{\i}ne holomorphe positive, implique qu'il existe
un voisinage connexe $V_z$ de $z$ tel que  $[\widetilde T]$
soit une $2$-cha\^{\i}ne holomorphe positive
 dans $V_z \times X$.
Mais alors, dans $(V_z \cap U_{max})\times X$, 
$[T]$ et $[\widetilde T]$
sont deux solutions positives au probl\`eme du bord pour $[\Gamma_f]$,
et donc elle sont \'egales dans $(V_z \cap U_{max})\times X$ 
d'apr\`es le lemme 3.21. D'apr\`es le lemme 3.8, on d\'eduit alors
que $[T]=[\widetilde T]$ dans tout $V_z \times X$. Mais
alors pour tout $x \in V_z \cap \widetilde V$, 
$[\widetilde S_x]=[\widetilde T, \pi,x]=[T,\pi,x]=[S_x]$
ce qui est impossible pas d\'efinition de $[\widetilde S_x]$
et donne la contradiction recherch\'ee.
\end{proof}

Montrons maintenant par l'absurde que $U_{max}=\CC$.
Supposons le contraire et notons  $G=\CC \backslash U_{max}$.
D'apr\`es le contr\^ole du volumme du lemme 3.7, $[T]$ est de 
masse finie et la masse des courants $[T,\Pi,z]$ est uniform\'ement born\'ee. 
L'ensemble $G$ ne peux donc \^etre non $0$-\'en\'erique 
(i.e non d\'enombrable).
En effet,  un ensemble ferm\'e et d\'enombrable admet toujours des points
isol\'es au voisinage desquels $[T]$ admettrait
une extension simple ferm\'ee ce qui contredirait la maximalit\'e de
$U_{max}$.
Supposons donc que $G$ est $0$-g\'en\'erique.
De la m\^eme mani\`ere que dans le lemme pr\'ec\'edent,
on montre qu'il existe alors un point $z \in \partial U_{max}$
et un voisinage $V_z$ de $z$ et un ouvert de Stein $W$
tel que le probl\`eme du bord
pour $[\Gamma]$ admette une solution positive $[\widetilde T]$ 
dans $W \cap (V_z \times X)$. Mais alors dans 
$(U_{max}\cap V_z)\times X$,
$[T]$ et $[\widetilde T]$ sont deux solutions positives 
dont l'union de leur support avec $\Gamma_f$ ne contient 
pas de sous-ensembles analytiques de dimension $2$,
on a donc, d'apr\`es le lemme 3.21,  
$[T]=[\widetilde T]$ dans $U_{max} \cap V_z$ 
et donc $[\Gamma_f]$ admet une solution au probl\`eme
du bord dans $(U_{max}\cup V_z)\times X$ 
ce qui contredit la maximalit\'e
de $U_{max}$ et montre que le probl\`eme du bord pour $[\Gamma_f]$
est r\'esoluble dans $\CC \times X$.
La projection de $[T]$ sur $X$ 
donne alors une solution au probl\`eme  du bord pour $[M]$ dans $X$. 
D'apr\`es \cite{henkin2}, $M$ est donc plongeable dans
$\CC^n$ et les fonctions CR sur $M$ s\'eparent les points.
\end{proof}\par

\end{document}